\documentclass[11pt]{article}

\usepackage{amsmath}
\usepackage{amsthm}
\usepackage{amsfonts}
\usepackage{amssymb}
\usepackage{mathtools}

\usepackage{hyperref}
\usepackage[compress,sort]{cite}

\usepackage{enumitem}

\usepackage{cases}
\usepackage{xcolor}

\newcommand{\dz}{\partial_z}
\newcommand{\dt}{\partial_t}
\newcommand{\dvh}{\mathrm{div}_h\,}

\newcommand{\nablah}{\nabla_h}
\newcommand{\Deltah}{\Delta_h}

\newcommand{\norm}[2]{\Arrowvert #1 \Arrowvert_{#2}}

\newcommand{\mrm}[1]{\mathrm{#1}}

\theoremstyle{plain}
\newtheorem{thm}{Theorem}[section]

\theoremstyle{definition}

\theoremstyle{remark}
\newtheorem{remark}{Remark}

\numberwithin{equation}{section}
\allowdisplaybreaks
\title{On the boundary layer arising from fast internal waves dynamics}

\author{
    Rupert Klein\footnote{
    {Institut f\"ur Mathematik}, {Freie Universit\"at Berlin}, {{Berlin}, {14195},  
    {Germany}}, {rupert.klein@math.fu-berlin.de}.
    }
    \qquad
    Xin Liu\footnote{
    {Department of Mathematics}, {Texas A\&M University}, 
    {College Station}, {77843},
    {Texas, U.S.},
    {xliu23@tamu.edu},
    {stleonliu@gmail.com}.}
    }

\begin{document}

\maketitle
% \tableofcontents

\begin{abstract}
    In this paper, we investigate the boundary layer arising from the fast internal waves in the Boussinesq equations with the Brunt-Vais\"al\"a frequency of order $ \mathcal O(1/\varepsilon) $. For the homogeneous-in-height stratification, previous work by \emph{Desjardins, Lannes, Saut, 3(1):153--192, Water Waves, 2021} establishes uniform-in-$\epsilon$ estimates locally in time, with additional constraints on the boundary data initially, which essentially restricts the dynamics in the spatially periodic domain. Removing such constraints, our goal is to investigate the general near-boundary behavior. We observe that the fast internal waves will give rise to large growth of the spatial derivatives in the normal direction of the solutions in the vicinity of the boundary. To capture this phenomenon, we introduce an inviscid boundary layer using a natural scaling. In addition, we investigate the well-posedness of such a boundary layer system in the space of analytic functions. 

	
\end{abstract}

\tableofcontents

\section{Internal waves and the Boussinesq equations}\label{sec:equations}

The propagation of internal waves in the ocean and atmosphere is affected by the stratification, which itself is the consequence of salinity, temperature, and moisture gradients (\cite{feliksLowFrequencyVariabilityMidlatitude2004, flierlModelsVerticalStructure1978, Klein2010}). There has been important progress in theoretical and experimental studies of internal waves. See, for instance, \cite{benjaminInternalWavesPermanent1967, davisSolitaryInternalWaves1967, longBoussinesqApproximationIts1965, onoAlgebraicSolitaryWaves1975, yihExactSolutionsSteady1960}. We refer to 
\cite{helfrichLongNonlinearInternal2006} for a review. 

Using {two-layer} flow models, Barros and Choi \cite{barrosInhibitingShearInstability2009, barrosRegularizingStronglyNonlinear2013} studied the free surface separating two flow layers, and showed the onset of Kelvin-Helmholtz instability for large internal waves. See also \cite{ebinIllposednessRayleightaylorHelmholtz1988}, where the Rayleigh-Taylor instability was also discussed. A regularized model was proposed by Barros and Choi as well. In a systematic way, Bona, Lannes, and Saut \cite{bonaAsymptoticModelsInternal2008} derived some new models with non-local operators for the {two-layer} flow interface. Bresch and Renardy \cite{breschWellposednessTwolayerShallowwater2011} studied the well-posedness of a model of two-layer flows, extending the result of Guyenne, Lannes, and Saut \cite{guyenneWellposednessCauchyProblem2010}. For more results, we refer readers to 
\cite{craigHamiltonianLongwaveExpansions2005,iguchiTwophaseFreeBoundary1997, jamesInternalTravellingWaves2001, kamotski2DRayleighTaylor2005, lannesWaterWavesProblem2013, lannesStabilityCriterionTwoFluid2013}.
For the investigation of compressible two-layer flows, we refer to \cite{Guo2011b, Jang2016, Wang2016, Wang2014, Wang2012}.

The study of continuously strongly stratified models, on the other hand, is quite limited. The local well-posedness of the Cauchy problem for density dependent Euler flows was investigated by Danchin \cite{danchinWellposednessIncompressibleDensitydependent2010}. For the two dimensional Boussinesq equations, the local well-posedness theory was developed by {Chae and Nam} \cite{chaeLocalExistenceBlowup1997}. In \cite{inciWellposednessInviscid2D2018}, Inci showed that the solution map of the inviscid Boussinesq equations is nowhere uniformly continuous. Danchin and Paicu \cite{danchinGlobalWellPosednessIssues2009} studied global solutions with Yudovich's type data with heat conductivity. See also \cite{lariosGlobalWellposedness2D2010, Larios2013} for flows with anisotropic viscosity, and \cite{kukavicaLongTimeBehavior2020} for a recent study of viscous flows. In addition, the shear-buoyancy instability has recently been investigated by Bedrossian, Bianchini, Coti Zelati, and Dolce \cite{bedrossianNonlinearInviscidDamping2021}. In particular, \cite{bedrossianNonlinearInviscidDamping2021}, together {with} \cite{Bedrossian2015}, provides evidence of an instability caused by internal gravity waves. All these studies are in the setting with a fixed Brunt-Vais\"al\"a frequency of order unity, i.e., with slow internal waves. In \cite{desjardinsNormalModeDecomposition2021}, {Desjardin et al.\ study flows with nonhomogeneous and homogeneous stratifications, and} the fast internal waves are investigated in terms of normal modes. See \cite{Charve-2020,Charve-2020-2,Charve-2022} for the study of dispersion of the quasi-geostrophic approximation in $ \mathbb R^3 $ and \cite{bresch-gerard-varet-grenier-planetarygeostrophic} for a study of the planetary geostrophic approximation.

On the other hand, the Boussinesq/pseudo-incompressible/anelastic models in general provide a good approximation to the full compressible flow in the low Mach number regime, even with fast internal waves. See \cite{Feireisl2007d, Feireisl2009a, Klein2005, Klein2009, Klein2010}. In particular, our previous work \cite{breschSoundproofModelAcoustic2022} shows that for a simplified, ideal, and weakly stratified system, the fast acoustic waves can be filtered out for low Mach number flows, and the incompressible model with fast internal waves provides a solid approximation to the original full compressible model/flow. Our result essentially considers these systems of equations on spatially periodic domains, including, in particular, periodicity in the vertical direction, to avoid handling the boundary behavior of the flow. To go beyond this, a better understanding of the dynamics of the flows near rigid boundaries is necessary, even for the incompressible model. This is the main goal and motivation of the present work. 

To be more precise, consider in $\Omega := \mathbb T^2 \times (0,1) $ the Boussinesq equations with fast internal waves, i.e., %\rupert{[There are changes of sign in (1.1b,c)!]}
\begin{subequations}\label{sys:boussinesq}
    \begin{align}
        \dt v + v \cdot \nablah v + w \dz v && + \dfrac{\nablah p}{\varepsilon} && && = 0, \label{eq:v} \\
        \dt w + v \cdot \nablah w + w \dz w && + \dfrac{\dz p}{\varepsilon} && {-} \dfrac{\theta}{\varepsilon} && = 0, \label{eq:w} \\
        \dt \theta + v \cdot \nablah \theta + w \dz \theta && && {+} \dfrac{w}{\varepsilon} && = 0, \label{eq:theta} \\
        \dvh v + \dz w && && && = 0. \label{eq:dv-free}
    \end{align}
In addition, on the boundary $ \lbrace z = 0,1 \rbrace $, we impose the condition of impermeability, i.e.,
\begin{equation}\label{bc:w}
    w\vert_{z=0,1} = 0.
\end{equation}
\end{subequations}
Here, $ \frac{1}{\varepsilon}$, for $\varepsilon \ll 1$, is the Brunt-Vais\"al\"a frequency of the internal waves, which we assume to be constant throughout the present analysis. 

In \cite[Theorem 2]{desjardinsNormalModeDecomposition2021}, under the assumption that
\begin{equation}\label{hp:lannes}
    \partial_z^{2k} \theta = \partial_z^{2k} w = \partial_z^{2k+1} v = 0, ~ k \in \mathbb N, \qquad \text{on} \qquad \lbrace z = 0,1 \rbrace, 
\end{equation}
a property that is invariant for the regular solutions, a uniform-in-$\varepsilon$ well-posedness result of system \eqref{sys:boussinesq} is established. Our goal of this paper to {investigate} system \eqref{sys:boussinesq} without the assumption \eqref{hp:lannes}. This will inform us about 
\begin{itemize}
    \item how assumption \eqref{hp:lannes} affects the solutions or, equivalently, what is the benefit of \eqref{hp:lannes} in the previous study;
    \item the general dynamics, in particular the boundary effect of solutions to system \eqref{sys:boussinesq}.
\end{itemize}

In particular, we will investigate the dynamics of solutions near the boundary, and introduce a boundary layer expansion. Our main theorem can be stated informally as follows:
\begin{thm}\label{thm:main-theorem-informal}
	As $ \varepsilon \rightarrow 0^+ $, there exists a boundary layer dynamics of system \eqref{sys:boussinesq} for general initial data without the assumption of \eqref{hp:lannes} {on the slow time scale, \emph{i.e.}, for $t = O(1)$}. Under the requirement that the boundary layer flow has no feedback to the internal flow, the boundary layer system is well-posed in the space of analytic functions. 
\end{thm}
\noindent
The main theorem will be stated in section \ref{sec:functional_space_bl} in full {detail (thm.~\ref{thm:b-l-detailed})}. {A multiple times scale analysis that simultaneously addresses the fast internal waves and the slow evolution of the boundary layer structure is deferred to future work.}

\bigskip

The rest of this paper is organized as follows. In section \ref{sec:linear_analysis}, by studying a linear system, we point out the `noise' arising from the boundary without the assumption \eqref{hp:lannes}, which is the motivation for introducing the boundary layer. In particular, the linear boundary layer is discussed in section \ref{subsec:linear-boundary-layer}. We refer readers to some physical remark about this boundary layer in section \ref{sec:physical-remarks}. In section \ref{sec:boundary_layer_nonlinear}, we introduce the nonlinear boundary layer expansion. Section \ref{sec:well-posedness-bl} is devoted to investigating the well-posedness theory of the boundary layer. A linear system is introduced in section \ref{sec:iota-approximation}, by which the problem of possible overdetermination is overcome. A contraction mapping argument is established in section \ref{sec:mapping}, leading to the well-posedness of the boundary layer equations. 

We emphasize that the boundary layer we introduce and construct is really confined near the boundary, and has no feedback to the internal flows. This, on the one hand, brings the possible overdetermination in our analysis (see section \ref{sec:well-posedness-bl}), and on the other hand, might not enjoy the most physical generality. The {latter} problem, and in particular a study of multiple-scales interactions between fast internal waves in the bulk of the atmosphere with slow boundary layer dynamics near the ground, will be left {for future studies}. 

\section{Motivation and linear analysis}\label{sec:linear_analysis}

\subsection{Linear system and motivation}

Consider firstly the following linear system: %\rupert{[$\exists$ changes of sign in (2.1b,c)!]}
\begin{subequations}\label{sys:linear}
    \begin{align}
        \dt v_l && + \dfrac{\nablah p_l}{\varepsilon} && &&  = 0, \label{eq:v-l} \\
        \dt w_l && + \dfrac{\dz p_l}{\varepsilon} && {-} \dfrac{\theta_l}{\varepsilon} && = 0, \label{eq:w-l} \\
        \dt \theta_l && && {+} \dfrac{w_l}{\varepsilon} && = 0, \label{eq:theta-l} \\
        \dvh v_l + \dz w_l && && && = 0, \label{eq:dv-free-l}\\
        w_l\vert_{z=0,1} && && && = 0. \label{bc:w-l}
    \end{align}
\end{subequations}

Assume for the moment that solutions exist and are smooth enough so that the following calculation makes sense. We would like to demonstrate the growth of $ (\partial_z^2 w_l)\vert|_{z=0,1} $ and $ (\partial_z^2 \theta_l)\vert|_{z=0,1} $. Indeed, taking $ \dvh $ of \eqref{eq:v-l} and $ \dz $ of \eqref{eq:w-l} and summing up the resultants, thanks to \eqref{eq:dv-free-l}, lead to
\begin{equation}\label{eq:p-l}
    (\Deltah + \partial_{zz}) p_l - \dz \theta_l = 0.
\end{equation}
Meanwhile, taking $ \partial_z^2 $ of \eqref{eq:w-l} yields
\begin{equation}\label{eq:l-001}
    \dt \partial_{z}^2 w_l + \dfrac{1}{\varepsilon} \partial_z(\partial_{zz} p_l - \partial_z \theta_l) = \dt \partial_{zz} w_l - \dfrac{1}{\varepsilon} \partial_z \Deltah p_l = 0.
\end{equation}
After replacing $ \dz p_l $ using \eqref{eq:w-l} in the second equality of \eqref{eq:l-001}, we arrive at 
\begin{equation}\label{eq:w-l-dzz}
    \dt \partial_{z}^2 w_l + \dt \Deltah w_l - \dfrac{1}{\varepsilon} \Deltah \theta_l = 0.
\end{equation}
Taking the trace of \eqref{eq:w-l-dzz} to the boundary $ \lbrace z = 0, 1 \rbrace $, thanks to \eqref{bc:w-l}, leads to
\begin{equation}\label{eq:w-l-dzz-trace}
    \dt (\partial_{z}^2 w_l)\vert_{z=0,1} - \dfrac{1}{\varepsilon} (\Deltah \theta_l)\vert_{z=0,1} = 0.
\end{equation}
To calculate $ (\Deltah \theta_l)\vert_{z=0,1} $, notice that \eqref{eq:theta-l} and \eqref{bc:w-l} imply that 
\begin{equation}
    (\Deltah \theta_l)\vert_{z=0,1}(t) \equiv (\Deltah \theta_l)\vert_{z=0,1}(t=0) = (\Deltah \theta_{l,\mrm{in}})\vert_{z=0,1}. 
\end{equation}
Hereafter, $ (\cdot)_\mrm{in} $ denotes the initial data (value at $ t= 0 $) of the corresponding quantity. 
Therefore, from \eqref{eq:w-l-dzz-trace}, one can evaluate that
\begin{equation}\label{ev:w-l-dzz-trace}
    (\partial_{z}^2 w_l)\vert_{z=0,1}(t) = (\partial_{z}^2 w_{l,\mrm{in}})\vert_{z=0,1} + \dfrac{t}{\varepsilon} (\Deltah \theta_{l,\mrm{in}})\vert_{z=0,1}.
\end{equation}

On the other hand, applying $ \partial_{z}^2 $ to \eqref{eq:theta-l} directly implies
\begin{equation}
    \dt \partial_{z}^2 \theta_l = -\dfrac{1}{\varepsilon} \partial_{z}^2 w_l,
\end{equation}
from which one can calculate that 
\begin{equation}\label{ev:theta-l-dzz-trace}
    (\partial_{z}^2 \theta_l)\vert_{z=0,1}(t) = (\partial_{z}^2 \theta_{l,\mrm{in}})\vert_{z=0,1} - \dfrac{t}{\varepsilon}(\partial_{z}^2 w_{l,\mrm{in}})\vert_{z=0,1} - \dfrac{t^2}{2\varepsilon^2} (\Deltah \theta_{l,\mrm{in}})\vert_{z=0,1}.
\end{equation}

\bigskip

From $ (\partial_z^2 w_l)\vert|_{z=0,1} $ and $ (\partial_z^2 \theta_l)\vert|_{z=0,1} $ as in 
\eqref{ev:w-l-dzz-trace} and \eqref{ev:theta-l-dzz-trace}, one observes that, near the boundary, the following scaling arises after initial time,
\begin{equation}\label{lw:l-scaling-01}
    \partial_z^2 w_l \simeq \dfrac{1}{\varepsilon} \quad \text{and} \quad \partial_z^2 \theta_l \simeq \dfrac{1}{\varepsilon^2}.
\end{equation}

As a remark, similar growth of $ (\dz v_l)\vert_{z=0,1} $ can be observed. Indeed, taking the $ z $-derivative of \eqref{eq:v-l} and subtracting the result with the horizontal gradient of \eqref{eq:w-l} lead to 
\begin{equation}
    \dt (\dz v_l - \nablah w_l) = - \dfrac{\nablah \theta_l}{\varepsilon}.
\end{equation}
Taking the trace of the equation above to the boundary yields
\begin{equation}
    \dt (\dz v_l)\vert_{z=0,1} = - \dfrac{(\nablah \theta_l)\vert_{z=0,1}}{\varepsilon},
\end{equation}
which implies that after initial time,
\begin{equation}
    (\dz v_l)\vert_{z=0,1} \simeq \dfrac{1}{\varepsilon},
\end{equation}
thanks to the fact that $ \dt (\nablah \theta_l)\vert_{z=0,1} = 0 $ (or equivalently, $ (\nablah \theta_l)\vert_{z=0,1}(t) \equiv (\nablah \theta_{l,\mrm{in}})\vert_{z=0,1} $ ).

\bigskip

Therefore, in order to capture the profiles of solutions near the boundary $ \lbrace z = 0 \rbrace $, it is {natural} to consider the following ansatz,
\begin{equation}\label{lw:l-scaling-02}
    w_l = \varepsilon W_l,\quad z = \varepsilon \eta,
\end{equation}
with $ W_l = \mathcal O(1) $. It is easy to verify that, in view of \eqref{lw:l-scaling-01}, that
\begin{equation}\label{lw:l-scaling-03}
    \partial_\eta^2 W_l = \varepsilon \dz^2 w_l = \mathcal O(1) \quad \text{and} \quad \partial_\eta^2 \theta_l = \varepsilon^2 \dz^2 \theta_l = \mathcal O(1).
\end{equation}
Moreover, let 
\begin{equation}\label{def:pressure-bl-001}
    q_l := \dfrac{1}{\varepsilon} p_l.
\end{equation}
Then, in terms of $ (v_l,W_l,\theta_l,q_l) $ and $ (x,y,\eta,t) $, one can derive from system \eqref{sys:linear} that
\begin{gather*}
    \dt v_l +  \nablah q_l = 0, \qquad 
    \varepsilon^2 \dt W_l + \partial_\eta q_l {-} \theta_l = 0, \\
    \dt \theta_l  {+} W_l = 0, \\
    \dvh v_l + \partial_\eta W_l = 0, \qquad
    W_l\vert_{z=0} = 0.
\end{gather*}

\subsection{Linear boundary layer equations}\label{subsec:linear-boundary-layer}

Motivated by the structure above, we propose the following boundary layer equation near $ \lbrace z = 0 \rbrace $:
\begin{subequations}\label{sys:bl-l-01}
    \begin{align}
        \dt v_{b,l} +  \nablah q_{b,l} = 0, \label{eq:v-l-bl-01} \\
        \partial_\eta q_{b,l} {-} \theta_{b,l} = 0, \label{eq:w-l-bl-01} \\
        \dt \theta_{b,l} {+} W_{b,l} = 0, \label{eq:theta-l-bl-01} \\
        \dvh v_{b,l} + \partial_\eta W_{b,l} = 0, \label{eq:dv-free-l-bl-01} \\
        W_{b,l}\vert_{\eta=0} = 0. \label{bc:l-bl-01}
    \end{align}
In addition, the matching condition at $ \eta = \infty $ is 
\begin{equation}\label{bc:infty-l-bl-01}
    v_{b,l}, \theta_{b,l} \rightarrow 0, \qquad \text{as} \quad \eta \rightarrow \infty.
\end{equation}
\end{subequations}
Furthermore, let
\begin{equation}\label{def:pressure-bl-001-2}
    \pi_{b,l}(x,y):= q_{b,l}(x,y,\eta) - \int_0^\eta \theta_{b,l} (x,y,s) \,ds,
\end{equation}
{which is} independent of the $ \eta $-variable, thanks to \eqref{eq:w-l-bl-01}. Then \eqref{eq:v-l-bl-01} can be written as 
\begin{equation*}\tag{\ref{eq:v-l-bl-01}'}\label{eq:001}
    \dt v_{b,l} + \nablah \pi_{b,l} + \int_0^\eta \nablah \theta_{b,l} (s) \,ds = 0,
\end{equation*}
with
\begin{equation*}\tag{\ref{eq:w-l-bl-01}'}\label{eq:001-1}
    \partial_\eta \pi_{b,l} = 0.
\end{equation*}

\bigskip 
\begin{remark}
A few remarks about the matching condition \eqref{bc:infty-l-bl-01}. Suppose, in general, that the far fields for $ v_{b,l} $ and $ \theta_{b,l}$ at $ \eta = \infty $ are given by
\begin{equation}\label{rm:ff-001}
    \lim_{\eta\rightarrow \infty}v_{b,l}(x,\eta,t) = v_{b,l,\infty}(x,t) \quad \text{and}\quad \lim_{\eta\rightarrow \infty}\theta_{b,l}(x,\eta,t) = \theta_{b,l,\infty}(x,t).  
\end{equation} 
Similarly, let the far fields of $ W_{b,l} $ and $ q_{b,l} $ be 
\begin{equation}\label{rm:ff-002}
    \lim_{\eta\rightarrow \infty}W_{b,l}(x,\eta,t) = W_{b,l,\infty}(x,t) \quad \text{and}\quad \lim_{\eta\rightarrow \infty}q_{b,l}(x,\eta,t) = q_{b,l,\infty}(x,t).  
\end{equation}

\smallskip
Then, \eqref{eq:w-l-bl-01} implies that
\begin{equation}\label{rm:ff-005}
    q_{b,l}(x,\eta,t) = - \int_\eta^\infty \theta_{b,l}(x,\iota,t) \,d\iota + q_{b,l,\infty}(x,t).
\end{equation} 
Therefore, in order to have a meaningful integration, one requires that
\begin{equation}\label{rm:ff-006}
    \theta_{b,l,\infty}(x,t) = \lim_{\eta\rightarrow\infty} \theta_{b,l}(x,\eta,t) = 0.
\end{equation}
On the other hand, thanks to \eqref{eq:dv-free-l-bl-01} and \eqref{bc:l-bl-01}, one has
\begin{equation}\label{rm:ff-003}
    W_{b,l,\infty}(x,t) = -\int_0^\infty \dvh v_{b,l}(x,\iota,t)\,d\iota,
\end{equation}
which, as before, requires that
\begin{equation}\label{rm:ff-004}
    \dvh v_{b,l,\infty}(x,t) = \lim_{\eta\rightarrow \infty }\dvh v_{b,l}(x,\eta,t) = 0.
\end{equation}
Meanwhile, \eqref{eq:v-l-bl-01} implies that
\begin{equation}\label{rm:ff-005}
    \dt v_{b,l,\infty} + \nablah q_{b,l,\infty} = 0.
\end{equation}
Consequently, $ v_{b,l}' := v_{b,l} - v_{b,l,\infty} ,  W_{b,l}' := W_{b,l} - W_{b,l,\infty},  \theta_{b,l}' := \theta_{b,l} - \theta_{b,l,\infty} $, and $ q_{b,l}' := q_{b,l} - q_{b,l,\infty} $ satisfy the same system of equations as $ v_{b,l}, W_{b,l}, \theta_{b,l}$, and $ q_{b,l} $, and satisfying the far field conditions \eqref{bc:infty-l-bl-01}. 

\end{remark}

Notice that, the matching condition of $ \theta_{b,l} $ at $ \eta \rightarrow \infty $ implicitly implies that 
\begin{equation}\label{bc:infty-l-bl-02}
    \lim_{\eta \rightarrow \infty }W_{b,l} = 0,
\end{equation}
thanks to \eqref{eq:theta-l-bl-01}. Therefore, the well-posedness of system \eqref{sys:bl-l-01}, even {though} it is linear, is not clear. This is due to the possible overdetermination of the far field condition \eqref{bc:infty-l-bl-01} (and \eqref{bc:infty-l-bl-02}). 

To resolve this problem, at the linear level, one should consider the following $ \iota $-approximating system:
\begin{subequations}\label{sys:bl-l-01-iota}
\begin{align}
    \dt v_{b,l,\iota} +  \nablah \pi_{b,l,\iota} + \int_0^\eta \nablah \theta_{b,l,\iota}(s)\,ds = 0, \label{eq:v-l-bl-01-iota} \\
    \partial_\eta \pi_{b,l,\iota} = 0, \label{eq:w-l-bl-01-iota} \\
    \dt \theta_{b,l,\iota} {+} W_{b,l,\iota} = 0, \label{eq:theta-l-bl-01-iota} \\
    \dvh v_{b,l,\iota} + \partial_\eta W_{b,l,\iota} = 0, \label{eq:dv-free-l-bl-01-iota} \\
    W_{b,l,\iota}\vert_{\eta=0,\iota} = 0, \label{bc:l-bl-01-iota} \\
    \text{in} ~ \Omega_\iota := \mathbb T^2 \times (0,\iota). \label{def:app-domain}
\end{align}
\end{subequations}
In particular, the approximating system \eqref{sys:bl-l-01-iota} replaces the domain of system \eqref{sys:bl-l-01} with a compact channel $ \Omega_\iota $ and the far field condition \eqref{bc:infty-l-bl-01} with $ W_{b,l,\iota}\vert_{\eta=\iota} = 0 $. The well-posedness theory of system \eqref{sys:bl-l-01-iota} is now classical. 
One can obtain proper uniform-in-$\iota$ estimates in some Sobolev space, $ H^k $ for instance.
Then after sending $ \iota \rightarrow \infty $, one can obtain the solution to system \eqref{sys:bl-l-01}. We remark that the far field condition \eqref{bc:infty-l-bl-01} is recovered thanks to the decay-in-$\eta$ property of $ H^k $ functions. 

\bigskip

Nevertheless, we assume for the moment that system \eqref{sys:bl-l-01} is well-posed. Then sending $ \eta \rightarrow \infty $ in \eqref{eq:001} and substituting \eqref{bc:infty-l-bl-01} yields that 
\begin{equation}\label{eq:far-field}
    \nablah \pi_{b,l} = - \int_0^\infty \nablah \theta_{b,l}(s)\,ds.
\end{equation}
Therefore, \eqref{eq:001} can be reduced to 
\begin{equation*}{\tag{\ref{eq:v-l-bl-01}'}}\label{eq:v-l-bl-01-02}
    \dt v_{b,l} - \nablah \int_\eta^\infty \theta_{b,l}(s)\,ds = 0.
\end{equation*}
Meanwhile, thanks to \eqref{eq:dv-free-l-bl-01} and \eqref{bc:l-bl-01}, \eqref{eq:theta-l-bl-01} can be written as
\begin{equation*}{\tag{\ref{eq:theta-l-bl-01}'}}\label{eq:theta-l-bl-01-02}
    \dt \theta_{b,l} - \int_0^\eta \dvh v_{b,l}(s)\,ds = 0.
\end{equation*}

\bigskip

It is easy to obtain the $ L^2 $ estimates of horizontal and time derivatives, i.e., for any $ k \in \mathbb N $, $ \sum_{\partial \in \lbrace \partial_x, \partial_y, \dt \rbrace }\norm{\partial^k v_{b,l}, \partial^k \theta_{b,l}}{L^2(\mathbb T^2\times (0,\infty))} \leq C_k $, for some $ C_k \in (0,\infty) $, depending only on the initial data. On the other hand, after applying $ \partial_\eta^k $ to \eqref{eq:v-l-bl-01-02} and \eqref{eq:theta-l-bl-01-02} and taking the $L^2$-inner product of the resultants with $ 2 \partial_\eta^k v_{b,l} $ and $ 2 \partial_\eta^k \theta_{b,l} $, respectively, summing up the estimates yields
\begin{equation}\label{est:l-001}
    \begin{gathered}
        \dfrac{d}{dt} \norm{\partial_\eta^k v_{b,l}, \partial_\eta^k \theta_{b,l}}{L^2(\mathbb T^2\times (0,\infty))}^2 \\
        = \int ( - \nablah \partial_\eta^{k-1} \theta_{b,l} \cdot \partial_\eta^k v_{b,l} + \dvh \partial_\eta^{k-1} v_{b,l} \partial_\eta^k \theta_{b,l}) \,dxdyd\eta \\
        =  \int \partial_\eta (\partial_\eta^{k-1}\theta_{b,l} \dvh \partial_\eta^{k-1} v_{b,l}) \,dxdyd\eta \\
        = - \int_{\mathbb{T}^2} (\partial_\eta^{k-1}\theta_{b,l} \dvh \partial_\eta^{k-1} v_{b,l})\vert_{\eta=0} \,dxdy \\
        \leq \norm{\partial_\eta^{k-1}\theta_{b,l}\vert_{\eta=0}}{H^{1/2}(\mathbb T^2)}\norm{\partial_\eta^{k-1}v_{b,l}\vert_{\eta=0}}{H^{1/2}(\mathbb T^2)} \\
        \lesssim \norm{\theta_{b,l}}{H^k(\mathbb T^2\times (0,\infty))}\norm{v_{b,l}}{H^k(\mathbb T^2\times (0,\infty))}.
    \end{gathered}
\end{equation}
All other mixed-derivatives can be estimated in a similar way. The estimates of $ W_{b,l} $ follows from \eqref{eq:theta-l-bl-01} and \eqref{eq:dv-free-l-bl-01}. In conclusion, one can establish, for any $ k \in \mathbb N $,
\begin{equation}\label{est:l-002}
    \begin{gathered}
    \sum_{ n+m = k} \norm{\dt^m v_{b,l}(t), \dt^m \theta_{b,l}(t)}{H^n(\mathbb T^2\times (0,\infty))} \\
    + \sum_{ n+m = k-1} \norm{\dt^m W_{b,l}(t), \dt^m \partial_\eta W_{b,l}(t)}{H^{n}(\mathbb T^2\times (0,\infty))}\\
    \leq e^{C_k t} \sum_{n+m=k} \norm{\dt^m v_{b,l,\mrm{in}}, \dt^m \theta_{b,l,\mrm{in}}}{H^n(\mathbb T^2\times (0,\infty))},
    \end{gathered}
\end{equation}
for some generic constant $ C_k \in (0,\infty) $, depending only on $ k $. Notice that, similarly to the Prandtl equations and the Hydrostatic Euler equations, \eqref{est:l-002} implies a loss of regularity of $ W_{b,l} $ in the horizontal derivatives. This indicates that one should consider analytic initial data in the corresponding nonlinear problem.

\subsection{Some physical remarks}\label{sec:physical-remarks}

The flow regime considered here corresponds to strongly stably stratified flow with stratification strength of the same order of magnitude throughout the troposphere. This is a typical situation found during night time under clear skies in the middle latitudes or even during the day in the polar regions. In the bulk of the troposphere, such stratification supports fast internal waves as discussed in the context of \eqref{sys:boussinesq} above. In contrast, in the boundary layer it gives rise to slow gravity currents as described by the present boundary layer analysis. 

A fully consistent solution of the bulk and boundary layer flows on the longer of the two time scales will require a systematic multiple time scale analysis, which we defer to future work. Nevertheless, we learn already from the present boundary layer equations that local heating from the ground -- corresponding to non-homogeneous $\theta|_{z=0}$ -- will not necessarily generate fast internal waves, but may merely induce slow near-ground motions if imposed on sufficiently slow time scales. 

Note also that the present boundary layer regime is decidedly different from the recent triple-deck boundary layer extension of quasi-geostrophic / Ekman layer theory in \cite{KleinEtAl2021}. The latter study addresses the coupling of quasigeostrophic flow to boundary layers with arbitrary stratification (including neutral) as they arise, e.g., in the morning at mid-latitudes when thermal convection tends to dissolves stable stratifications that may have developed over night. 

% =======================================================
% =======================================================
% =======================================================

\section{Boundary layer expansion: Formulation of nonlinear equations}\label{sec:boundary_layer_nonlinear}

We return to system \eqref{sys:boussinesq}. In particular, we consider solutions to system \eqref{sys:boussinesq} that consist of two parts: the internal bulk part and the boundary layer near $ z = 0 $. That is, we assume the following ansatz, 
\begin{equation}\label{def:bulk-boundary-dcmp}
    (v, w, \theta) = (v_{p},w_p,\theta_p)(x,y,z,t) + (v_b,\varepsilon W_b,\theta_b)(x,y,\dfrac{z}{\varepsilon},t) + \mathrm{Err},
\end{equation}
with the following assumptions
\newcounter{assumption}\setcounter{assumption}{0}
\begin{enumerate}[label = {\bf H\arabic*)},ref = {\bf H\arabic*)}]
    \setcounter{enumi}{\value{assumption}}
    \item \label{hp:h1} $ (v_p, w_p,\theta_p) $ satisfies \eqref{hp:lannes};
    \item \label{hp:h2} $ \lim_{\eta\rightarrow +\infty} (v_b,\theta_b)(x,y,\eta,t) = 0 $; 
    \item \label{hp:h3} $ \mathrm{Err} \rightarrow 0 $ as $ \varepsilon \rightarrow 0^+ $. 
    \setcounter{assumption}{\value{enumi}}
\end{enumerate}
Meanwhile, we denote by 
\begin{equation}\label{def:bulk-boundary-dcmp-p}
    p = p_p(x,y,z,t) + \varepsilon p_b(x,y,\dfrac{z}{\varepsilon}, t).
\end{equation}
We remark that \ref{hp:h1}--\ref{hp:h3} are important assumptions to formally derive the boundary layer equations. To justify \eqref{def:bulk-boundary-dcmp}, one need to prove \ref{hp:h3}. However, for the moment, we leave such issues for later. Also, we will not present the boundary layer near $ z = 1 $, which has the same structure as the one near $ z = 0 $ and is omitted for the sake of presentation.

\bigskip

Then, the internal bulk system reads, in $ \Omega $,
\begin{subequations}\label{sys:bulk_eqs}
    \begin{align}
        \dt v_p + v_p \cdot \nablah v_p + w_p \dz v_p && + \dfrac{\nablah p_p}{\varepsilon} && && = 0, \label{eq:bulk-v} \\
        \dt w_p + v_p \cdot \nablah w_p + w_p \dz w_p && + \dfrac{\dz p_p}{\varepsilon} && {-} \dfrac{\theta_p}{\varepsilon} && = 0, \label{eq:bulk-w} \\
        \dt \theta_p + v_p \cdot \nablah \theta_p + w_p \dz \theta_p && && {+} \dfrac{w_p}{\varepsilon} && = 0, \label{eq:bulk-theta} \\
        \dvh v_p + \dz w_p && && && = 0, \label{eq:bulk-dv-free}
    \end{align}
    with
    \begin{equation}\label{bc:bulk}
        \partial_z^{2k} \theta_p = \partial_z^{2k} w_p = \partial_z^{2k+1} v_p = 0, ~ k \in \mathbb N, \qquad \text{on} \qquad \lbrace z = 0,1 \rbrace, 
    \end{equation}
\end{subequations}
We remind readers that the boundary condition of system \eqref{sys:bulk_eqs} is the same as in \eqref{bc:w}, and \eqref{bc:bulk} is invariant as proved in \cite{desjardinsNormalModeDecomposition2021}. 

Meanwhile, to obtain the boundary layer equations, after substituting \eqref{def:bulk-boundary-dcmp} and \eqref{def:bulk-boundary-dcmp-p} in system \eqref{sys:boussinesq}, one can write down the equations in the coordinate $ (x,y,\eta :=\frac{z}{\varepsilon},t)\in \mathbb T^2 \times [0,\infty) \times [0,\infty) $. That is,
\begin{subequations}\label{sys:boundary_eqs}
\begin{gather}
    \dt v_b + (v_p\vert_{z=0} + v_b) \cdot \nablah v_b + v_b \cdot \nablah v_p\vert_{z=0} + W_b \partial_\eta v_b + \nablah p_b = 0, \label{eq:boundary-v} \\
    % \varepsilon^2(\dt W_b + (v_p\vert_{z=0} + v_b) \cdot \nablah W_b  + W_b \partial_\eta W_b + W_b (\dz w_p)\vert_{z=0} ) +
    \partial_\eta p_b - \theta_b = 0, \label{eq:boundary-w}\\
    \dt \theta_b + (v_p\vert_{z=0} + v_b) \cdot \nablah \theta_b + W_b \partial_\eta \theta_b + W_b = 0, \label{eq:boundary-theta} \\
    \dvh v_b + \partial_\eta W_b = 0, \label{eq:boundary-dv-free}
\end{gather}
\end{subequations}
where the bulk variables are valued at the boundary, i.e., $ z = 0 $, the trace of system \eqref{sys:bulk_eqs} at $ z = 0 $ is used, and the $ \mathrm{Err} $ and $ \mathcal O(\varepsilon) $ terms are omitted.  

\bigskip

The main goal of this paper is to investigate the well-posedness of the boundary layer system \ref{sys:boundary_eqs}. This will be done in the next section. 

\section{Well-posedness theory of the boundary layer equations, i.e., \eqref{sys:boundary_eqs}}\label{sec:well-posedness-bl}

In this section, we will investigate the well-posedness of the boundary layer equations, i.e., \eqref{sys:boundary_eqs}, arising from the boundary layer expansion in section \ref{sec:boundary_layer_nonlinear}. For this purpose, throughout this section, we assume that $ v_p\vert_{z=0} $ is given, bounded in $ X^{s+1} $, $ s \geq 3 $, defined in \eqref{def:functional-s-iota}, below. 

To simplify the notation, by denoting
\begin{equation}
    v = v_b, \quad w = W_b, \quad p = p_b, \quad \theta = \theta_b, \quad \text{and} \quad V = v_p\vert_{z=0},
\end{equation}
we rewrite system \eqref{sys:boundary_eqs} as 
\begin{subequations}\label{sys:bl}
    \begin{gather}
        \dt v + (V + v) \cdot \nablah v + w \partial_\eta v + v \cdot \nablah V + \nablah p = 0, \label{eq:bl-v}\\
        \partial_\eta p - \theta = 0, \label{eq:bl-w}\\
        \dt \theta + (V + v) \cdot \nablah \theta + w \partial_\eta \theta + w = 0, \label{eq:bl-theta}\\
        \dvh v + \partial_\eta w = 0, \label{eq:bl-dv-free}
    \end{gather}
with 
\begin{equation}\label{bc:bl-v-theta}
    \lim_{\eta\rightarrow \infty}(v,\theta)(x,y,\eta,t) = 0 \quad \text{and} \quad w\vert_{\eta=0} = 0.
\end{equation}
\end{subequations}

\bigskip

Before we continue, we would like to point out the difficulty of studying the well-posedness of solutions to system \eqref{sys:bl}. First, the same difficulty as the linear system \eqref{sys:bl-l-01} as discussed in section \ref{subsec:linear-boundary-layer} is still present in the nonlinear system \eqref{sys:bl}. That is, the possible overdetermination of $ \lim_{\eta\rightarrow \infty}w = 0 $ due to \eqref{bc:bl-v-theta} and \eqref{eq:bl-theta}. As in system \eqref{sys:bl-l-01-iota}, this should be resolved by introducing an $\iota$-approximating, hydrostatic system in domain with finite depth. Then similar to \eqref{est:l-001} and \eqref{est:l-002}, this shall enable uniform-in-$ \iota $ $ H^s $ estimates for some $ s \in \mathbb Z $. However, due to the nonlinearity, in particular, the vertical transport terms $ w \partial_\eta v $ and $ w \partial_\eta \theta $, two problems may arise:
\begin{enumerate}
    \item Loss of horizontal regularity;
    \item Loss of boundedness of $ w $. 
\end{enumerate}
Both these problems are due to the fact that the vertical velocity $ w $ is recovered through integral in the $ \eta $-variable of the incompressible condition \eqref{eq:bl-dv-free}, i.e., $ w = - \int_0^\eta \dvh v(s)\,ds $. While the first problem can be resolved by considering analytic-in-$ \vec x_h = (x,y) $ regularity of the solutions, the second problem is usually overcome by introducing decay in $ \eta $, or equivalently, a $ \eta $-weighted norm of the solutions. One can observe that, from the linear structure \eqref{eq:v-l-bl-01-02} and \eqref{eq:theta-l-bl-01-02}, any algebraic decay in $ \eta $ cannot be balanced in the equations thanks to the non-local internal wave terms, i.e., $ \nablah \int_\eta^\infty \theta_{b,l}(s)\,ds $ and $ \int_0^\eta \dvh v_{b,l}(s)\,ds $. Therefore a natural strategy is to consider exponential decay in $ \eta $. However, one can at most obtain algebraic decay in $ \iota $ of the pressure $ \pi_{b,l,\iota} $ in the $\iota$-approximating system \eqref{sys:bl-l-01-iota} (similarly in the corresponding nonlinear system). This is thanks to \eqref{eq:v-l-bl-01-iota}, \eqref{eq:v-l-bl-01-iota}, \eqref{eq:dv-free-l-bl-01-iota}, and \eqref{bc:l-bl-01-iota}, that
\begin{equation}
    \Deltah \pi_{b,l,\iota} = - \dfrac{1}{\iota} \int_0^\iota \int_0^\eta \Deltah \theta_{b,l,\iota}(s)\,dsd\eta.
\end{equation}
This makes it almost impossible to obtain uniform-in-$\iota $ estimates of solutions with exponential decay in $ \eta $ to the $ \iota $-approximating system. 

To circumvent these problems, our strategy is the following:
\begin{enumerate}[leftmargin=1.5cm]
    \item[Step 1.] Introducing an $ \iota $-approximating, hydrostatic system with frozen transport speed, and obtaining the uniform-in-$\iota $ $ H^s $ estimates;
    \item[Step 2.] Sending $ \iota \rightarrow \infty $, and obtaining a solution to the system in $ \Omega^+ $ without the pressure term;
    \item[Step 3.] Using fixed point arguments in functional space with analyticity and exponential decay in $ \eta $ to obtain solutions to system \eqref{sys:bl}. 
\end{enumerate}
We would like to point out that 1. the same difficulties do not exist in the study of the Prandtl equations, since there is no possible overdetermination as in \eqref{bc:infty-l-bl-02} in the Prandtl equations; 2. thanks to the study of Euler equations in the Eulerian and Lagrangian coordinates in 
\cite{constantinContrastLagrangianEulerian2016}, even though the loss of regularity only appears in the horizontal variables, it is unlikely to obtain well-posedness in the anisotropic analytic space. 

\subsection{The $\iota$-approximating system with frozen transport speed and uniform-in-$\iota$ estimates}\label{sec:iota-approximation}

As mentioned before, we shall introduce the $\iota$-approximating system with frozen transport speed in this section. 
% Due to the non-local term and the possible overdetermination as mentioned in \eqref{bc:infty-l-bl-02}, we introduce an approximation of system \eqref{sys:bl} below, 
Similar to system \eqref{sys:bl-l-01-iota}, we consider the following system of equations: 
in $ \Omega_\iota\times \lbrace t \geq 0 \rbrace $ with $ \Omega_\iota := \mathbb T^2 \times (0,\iota) = \lbrace (\vec x_h,\eta) \vert \vec x_h = (x,y)\in \mathbb T^2, \eta \in (0,\iota) \rbrace $ and $ \iota > 0 $, 
\begin{subequations}\label{sys:bl-app}
    \begin{gather}
        \dt v_\iota + (V + v_\iota^\mrm o) \cdot \nablah v_\iota + w_\iota^\mrm o \partial_\eta v_\iota + v_\iota \cdot \nablah V + \nablah p_\iota = 0,\label{eq:bl-app-v} \\
        \partial_\eta p_\iota - \theta_\iota  = 0, \label{eq:bl-app-w} \\
        \dt \theta_\iota + (V + v_\iota^\mrm o) \cdot \nablah \theta_\iota + w_\iota^\mrm o \partial_\eta \theta_\iota + w_\iota = 0, \label{eq:bl-app-theta} \\
        \dvh v_\iota + \partial_\eta w_\iota = 0, \label{eq:bl-app-dv-free}
    \end{gather}
    with
    \begin{equation}\label{bc:bl-app}
        w_\iota\vert_{\eta = 0,\iota} = 0.
    \end{equation}
\end{subequations}
Here, $ (v_\iota^\mrm o, w_\iota^\mrm o ) $ is given and satisfies 
\begin{equation}\label{eq:transport}
    \dvh v_\iota^\mrm o + \partial_\eta w_\iota^\mrm o = 0 \qquad\text{and}\qquad w_\iota^\mrm o\vert_{\eta = 0,\iota} = 0.
\end{equation}
We remind readers that $ V = V(x,y,t) $ is the trace of the internal bulk flow and is assumed to be given.

\bigskip

Moreover, by writing
\begin{equation}\label{def:varpi-app}
    \pi_\iota(\vec x_h) := p_\iota(\vec x_h, \eta) - \int_0^\eta \theta_\iota(\vec x_h, s) \,ds,
\end{equation}
\eqref{eq:bl-app-v} and \eqref{eq:bl-app-w} can be written as
\begin{gather*}
    \begin{gathered}
        \dt v_\iota + (V + v_\iota^\mrm o) \cdot \nablah v_\iota + w_\iota^\mrm o \partial_\eta v_\iota + v_\iota \cdot \nablah V + \nablah \pi_\iota \\
        + \int_0^\eta \nablah \theta_\iota(s)\,ds = 0,
    \end{gathered}
    \tag{\ref{eq:bl-app-v}'} \label{eq:bl-app-v-1} \\
    \text{and} \qquad\qquad\qquad \partial_\eta \pi_\iota = 0. \tag{\ref{eq:bl-app-w}'} \label{eq:bl-app-w-1}
\end{gather*}
In addition, \eqref{bc:bl-app} and \eqref{eq:bl-app-dv-free} imply that
\begin{equation}\label{eq:bc-app-w-1}
    w_\iota(\eta) = - \int_0^\eta \dvh v_\iota(s)\,ds = \int_\eta^\iota \dvh v_\iota(s)\,ds,
\end{equation}
and
\begin{equation}\label{eq:bc-app-w-2}
    \int_0^\iota \dvh v_\iota(s)\,ds = 0.
\end{equation}

\bigskip

Next, we shall obtain the uniform-in-$ \iota $ estimates of solutions to system \eqref{sys:bl-app}. Let $ s \in \mathbb Z^+ $ be any positive integer, and consider the functional space
\begin{equation}\label{def:functional-s-iota}
    X^s = X^s_t:= \lbrace f = f(x,y,\eta,t) \vert f \in H^s(\Omega_\iota), ~ \dt f \in H^{s-1}(\Omega_\iota) \rbrace
\end{equation}
and the associated $ X^s $ norm 
\begin{equation}\label{def:X-norm-iota}
    \norm{f(t)}{X^s}:= \sqrt{ \norm{f(t)}{H^s(\Omega_\iota)}^2 + \norm{\dt f(t)}{H^{s-1}(\Omega_\iota)}^2 }.
\end{equation}

\bigskip

Let $ \partial_h \in \lbrace \partial_x, \partial_y \rbrace $ and $ \partial \in \lbrace \dt, \partial_x, \partial_y, \partial_\eta \rbrace $. For any multi-index $ \mathbf s = (s_0,s_1, s_2) $ with $ s_0 \in \lbrace 0,1 \rbrace $, $ s_1,s_2 \in \mathbb N $, and $ \vert \mathbf s \vert = s_0+ s_1 + s_2 \leq s $, denote by 
\begin{equation}\label{def:multi-index-derivative}
    \partial^\mathbf s := \partial_\eta^{s_2} \partial_{h}^{s_1} \dt^{s_0}.
\end{equation}
In addition, we also denote the $L^2$-inner product by 
\begin{equation}\label{def:inner-product-iota}
    \langle f,g  \rangle := \int_{\Omega_\iota}f \cdot g \,dxdyd\eta.
\end{equation}
Then after applying $ \partial^\mathbf s $ to \eqref{eq:bl-app-v-1} and \eqref{eq:bl-app-theta} and taking the $ L^2 $-inner product of the resultant equations with $ 2 \partial^\mathbf s v_\iota $ and $ 2 \partial^\mathbf s \theta_\iota $, thanks to \eqref{eq:transport}, one arrives at
\begin{equation}\label{est:001}
    \begin{aligned}
        & \dfrac{d}{dt} \norm{\partial^\mathbf s v_\iota, \partial^\mathbf s \theta_\iota}{L^2(\Omega_\iota)}^2 = \underbrace{2 \langle \partial^\mathbf s \pi_\iota, \dvh \partial^\mathbf s v_\iota \rangle}_{=:I_1} \\
        & \qquad \underbrace{- 2\langle \partial^\mathbf s \int_0^\eta \nablah  \theta_\iota(s)\,ds,  \partial^\mathbf s v_\iota \rangle - 2\langle \partial^\mathbf s w_\iota, \partial^\mathbf s \theta_\iota \rangle}_{=:I_2} \\
        & \underbrace{- \sum_{f\in\lbrace v_\iota, \theta_\iota\rbrace}2 \langle \partial^\mathbf s\lbrack (V + v_\iota^\mrm o) \cdot \nablah f \rbrack - (V + v_\iota^\mrm o) \cdot \nablah \partial^\mathbf s f, \partial^\mathbf s f \rangle}_{=:I_3} \\
        & \qquad \underbrace{- \sum_{f\in \lbrace v_\iota, \theta_\iota\rbrace } 2 \langle \partial^\mathbf s ( w_\iota^\mrm o \partial_\eta f) - w_\iota^\mrm o \partial_\eta \partial^\mathbf s f, \partial^\mathbf s f\rangle}_{=:I_4} \\
        & \qquad + \underbrace{\int \dvh V \vert \partial^\mathbf s v_\iota \vert^2 \,dxdyd\eta - 2 \langle \partial^\mathbf s(v_\iota \cdot \nablah V), \partial^\mathbf s v_\iota \rangle}_{=:I_5}.
    \end{aligned}
\end{equation}
Now we estimate the right hand side of \eqref{est:001}. For $ I_1 $, thanks to \eqref{eq:bl-app-w-1}, $ I_1 = 0 $ if $ s_2 > 0 $. Otherwise, substituting \eqref{eq:bc-app-w-2} yields $ I_1 =0 $. 

To calculate $ I_2 $, when $ s_2 = 0 $, applying integration by parts yields that 
\begin{equation}
    \begin{aligned}
        I_2 = &  2 \langle \int_0^\eta \partial^\mathbf s \theta_\iota(s)\,ds, \partial^\mathbf s \dvh v_\iota \rangle - 2 \langle \partial^\mathbf s w_\iota, \partial^\mathbf s \theta_\iota \rangle \\
        = & - 2 \langle \int_0^\eta \partial^\mathbf s \theta_\iota(s)\,ds, \partial^\mathbf s \partial_\eta w_\iota \rangle - 2 \langle \partial^\mathbf s w_\iota, \partial^\mathbf s \theta_\iota \rangle \\
        = &  2 \langle \partial^\mathbf s \theta_\iota, \partial^\mathbf s  w_\iota \rangle - 2 \langle \partial^\mathbf s w_\iota, \partial^\mathbf s \theta_\iota \rangle = 0,  
    \end{aligned}
\end{equation}
thanks to \eqref{eq:bl-app-dv-free} and \eqref{bc:bl-app}.
Otherwise, when $ s_2 > 0 $, one has that
\begin{equation}
    \begin{aligned}
        I_2 = & - 2 \langle \nablah \partial_\eta^{s_2-1}\partial_h^{s_1} \partial_t^{s_0} \theta_\iota, \partial^\mathbf s v_\iota \rangle - 2 \langle \partial_\eta^{s_2-1} \partial_h^{s_1} \partial_t^{s_0} \partial_\eta w_\iota, \partial^\mathbf s \theta_\iota \rangle \\
        = &  2 \langle \partial_\eta^{s_2-1}\partial_h^{s_1} \partial_t^{s_0} \theta_\iota, \partial^\mathbf s \dvh v_\iota \rangle + 2 \langle \partial_\eta^{s_2-1} \partial_h^{s_1} \partial_t^{s_0} \dvh v_\iota, \partial^\mathbf s \theta_\iota \rangle \\
        = & 2 \int_{\Omega_\iota} \partial_\eta (\partial_\eta^{s_2-1}\partial_h^{s_1} \partial_t^{s_0} \theta_\iota \cdot \partial_\eta^{s_2-1}\partial_h^{s_1} \partial_t^{s_0} \dvh v_\iota  ) \,dxdy\eta \\
        = & 2 \int_{\mathbb{T}^2}(\partial_\eta^{s_2-1}\partial_h^{s_1} \partial_t^{s_0} \theta_\iota \cdot \partial_\eta^{s_2-1}\partial_h^{s_1} \partial_t^{s_0} \dvh v_\iota)(\eta)\Big\vert^\iota_0 \,dxdy \\
        \lesssim & \sum_{\eta = 0,\iota}\norm{\partial_\eta^{s_2-1}\partial_h^{s_1} \partial_t^{s_0} \theta_\iota(\eta)}{H^{1/2}(\mathbb T^2)} 
        \norm{\partial_\eta^{s_2-1}\partial_h^{s_1} \partial_t^{s_0} v_\iota(\eta)}{H^{1/2}(\mathbb T^2)} \\
        \lesssim & \norm{\theta_\iota}{X^s} \norm{v_\iota}{X^s}.
    \end{aligned}
\end{equation}

Applying Leibniz's rule to $ I_3 $ and $ I_4 $ implies
\begin{equation}
    \begin{aligned}
        I_3 + I_4 \lesssim & \biggl\lbrace \sum_{f \in \lbrace v_\iota, \theta_\iota \rbrace} \norm{\partial V, \partial v_\iota^\mrm o, \partial w_\iota^\mrm o}{L^\infty(\Omega_\iota)} \norm{\partial^{s-1} \nablah f, \partial^{s-1} \partial_\eta f}{L^2(\Omega_\iota)} \\
        % & \qquad\qquad \times \norm{\partial^s f}{L^2(\Omega_\iota)} \\
        + & \sum_{\substack{f \in \lbrace v_\iota, \theta_\iota \rbrace \\ 2 \leq l \leq s-1}} \norm{\partial^l V, \partial^l v_\iota^\mrm o, \partial^l w_\iota^o}{L^4(\Omega_\iota)}\norm{\partial^{s-l}\nablah f, \partial^{s-l} \partial_\eta f}{L^4(\Omega_\iota)} \\
        % & \qquad\qquad \times \norm{\partial^s f}{L^2(\Omega_\iota)} \\
        + & \sum_{f \in \lbrace v_\iota, \theta_\iota \rbrace} \norm{\partial^s V, \partial^s v_\iota^\mrm o, \partial^s w_\iota^\mrm o}{L^2(\Omega_\iota)} \norm{\nablah f, \partial_\eta f}{L^\infty(\Omega_\iota)} \biggr\rbrace \\
        & \qquad \times \norm{\partial^s v_\iota, \partial^s \theta_\iota}{L^2(\Omega_\iota)} \lesssim \norm{V, v_\iota^\mrm o, w_\iota^\mrm o}{X^s} \norm{v_\iota, \theta_\iota}{X^s}^2,
    \end{aligned}
\end{equation}
provided $ s \geq 3 $.

Last but not least, directly calculation shows that
\begin{equation}
    \begin{aligned}
        I_5 \lesssim & \norm{\nablah V}{L^\infty(\Omega_\iota)} \norm{\partial^s v_\iota}{L^2(\Omega_\iota)}^2 \\
        + & \sum_{1\leq l \leq s-1} \norm{\partial^{s-l}v_\iota}{L^4(\Omega_\iota)} \norm{\nablah \partial^l V}{L^4(\Omega_\iota)} \norm{\partial^s v_\iota}{L^2(\Omega_\iota)} \\
        + & \norm{v_\iota}{L^\infty(\Omega_\iota)}\norm{\nablah \partial^s V}{L^2(\Omega_\iota)}\norm{\partial^s v_\iota}{L^2(\Omega_\iota)} \\
        \lesssim & \norm{V}{X^{s+1}}\norm{v_\iota}{X^s}^2,
    \end{aligned}
\end{equation}
provided $ s \geq 2 $. 

In summary, from \eqref{est:001}, one can conclude from the estimates above that, for any $ s \geq 3 $,
\begin{equation}\label{est:002}
    \dfrac{d}{dt} \norm{v_\iota,\theta_\iota}{X^s}^2 \leq C_s ( 1 + \norm{V}{X^{s+1}}+ \norm{v_\iota^\mrm o, w_\iota^\mrm o}{X^s} ) \times
    \norm{v_\iota,\theta_\iota}{X^s}^2,
\end{equation}
with some constant $ C_s \in (0,\infty) $ independent of $ \iota $. 

\bigskip

To finish this subsection, let $ \Omega^+ := \mathbb T^2 \times (0,\infty) $. For $ f $ defined in $ \Omega_\iota $, we trivially extend $ f $ to $ \Omega_+ $ by zero in $ \Omega^+/\Omega_\iota $. Using such notations, we choose $ v_\iota^\mrm o $ and $ w_\iota^\mrm o $ such that, 
\begin{equation}\label{def:app-v-w-o}
    \begin{aligned}
        v_\iota^\mrm o \rightarrow v^\mrm o && \text{in} && L^\infty(0,\infty;H^s(\Omega_+))\cap H^1(0,\infty;H^{s-1}(\Omega_+)), \\
        \text{and} \qquad w_\iota^\mrm o \rightarrow w^\mrm o && \text{in} && L^\infty(0,\infty;H^s(\Omega_+))\cap H^1(0,\infty;H^{s-1}(\Omega_+)),
    \end{aligned}
\end{equation}
for some $ v^\mrm o, w^\mrm o \in L^\infty(0,\infty;H^s(\Omega_+))\cap H^1(0,\infty;H^{s-1}(\Omega_+)) $, and
\begin{equation}\label{def:app-v-w-o-1}
    \sup_{0\leq t < \infty} \biggl( \norm{v_\iota^\mrm o(t), w_\iota^\mrm o(t)}{H^s(\Omega_+)} + \norm{\dt v_\iota^\mrm o(t), \dt w_\iota^\mrm o(t)}{H^{s-1}(\Omega_+)} \biggr) \leq C
\end{equation}
for some constant $ C \in (0,\infty) $, independent of $ \iota $. Then from \eqref{est:002}, one can conclude that
\begin{equation}\label{est:003}
    \sup_{0\leq t \leq T} \norm{v_\iota, \theta_\iota}{X^s}^2 \leq e^{C_s(1+\norm{V}{X^{s+1}} + C)T} \times \norm{v_{\iota,\mrm{in}},\theta_{\iota, \mrm{in}}}{X^s}^2,
\end{equation}
for any $ T\in (0,\infty) $ and $ C $ independent of $ \iota $, where $ v_{\iota,\mrm{in}} $ and $ \theta_{\iota,\mrm{in}} $ are the initial data for $ v_\iota $ and $ \theta_\iota $, respectively. Therefore, consider the initial data such that
\begin{equation}\label{est:004}
    \begin{aligned}
        && (v_{\iota,\mrm{in}},\theta_{\iota,\mrm{in}}) \rightarrow (v_{\mrm{in}},\theta_{\mrm{in}}) && \text{in} & \quad H^s(\Omega_+)\\
        \text{and} && (\dt v_{\iota, \mrm{in}}, \dt \theta_{\iota,\mrm{in}} ) \rightarrow (v_{\mrm{in},1},\theta_{\mrm{in},1}) && \text{in} & \quad H^{s-1}(\Omega_+),
    \end{aligned}
\end{equation}
where $ \dt v_{\iota, \mrm{in}} $ and $ \dt \theta_{\iota,\mrm{in}} $ are defined by the equations \eqref{eq:bl-app-v-1} and \eqref{eq:bl-app-theta}, respectively, for some $ v_\mrm{in},\theta_\mrm{in} \in H^s(\Omega_+) $ and $ v_{\mrm{in},1},\theta_{\mrm{in},1} \in H^{s-1}(\Omega_+) $. Then similar bound as in \eqref{def:app-v-w-o-1} holds for $ v_{\mrm{in}},\theta_{\mrm{in}}, v_{\mrm{in},1},\theta_{\mrm{in},1} $.  
In addition, $ v^\mrm o, w^\mrm o, v_{\iota,\mrm{in}} $ satisfy the compatibility condition 
\begin{equation}\label{est:005}
    \begin{gathered}
        \dvh v^\mrm o + \partial_\eta w^\mrm o = 0, \qquad \int_0^\infty \dvh v^\mrm{o}(\eta)\,d\eta = \int_0^\infty \dvh v_{\iota, \mrm{in}}(\eta)\,d\eta = 0,\\
        \text{and} \qquad w^\mrm o\vert_{\eta=0}=0.
    \end{gathered}
\end{equation}
Therefore, we obtain a sequence of solutions $ (v_\iota, \theta_\iota) $ to system \eqref{sys:bl-app} with uniform bounds in $ L_\mrm{loc}^\infty(0,\infty;H^s(\Omega_+))\cap H^1_\mrm{loc}(0,\infty;H^{s-1}(\Omega_+)) $. 

With such uniform bound, by choosing subsequences if necessary, we obtain a solution to the following system: in $ \Omega_+ \times\lbrace t\geq 0 \rbrace $, 
\begin{subequations}\label{sys:bl-app-1}
    \begin{gather}
        \begin{gathered}
            \dt v + (V + v^\mrm o)\cdot \nablah v + w^\mrm o \partial_\eta v + v \cdot \nablah V + \nablah \pi \\
            + \int_0^\eta \nablah \theta (s)\,ds = 0,\end{gathered}\label{eq:bl-app-1-v}\\
        \partial_\eta \pi = 0, \label{eq:bl-app-1-w} \\
        \dt \theta + (V + v^\mrm o) \cdot \nablah \theta + w^\mrm o \partial_\eta \theta + w = 0, \label{eq:bl-app-1-theta} \\
        \dvh v + \partial_\eta w = 0, \label{eq:bl-app-1-dv-free}
    \end{gather}
    with
    \begin{equation}\label{bc:bl-app-1}
        w\vert_{\eta = 0} = 0\qquad \text{and} \qquad \int_0^\infty \dvh v(\eta)\,d\eta = 0. 
    \end{equation}
\end{subequations}
Moreover, $ (v,\theta) $ satisfies, for any $ s \geq 3$,
\begin{equation}\label{est:006}
    \sup_{0 \leq t \leq T} \biggl( \norm{v(t), \theta(t)}{H^s(\Omega_+)} + \norm{\dt v(t), \dt \theta(t)}{H^{s-1}(\Omega_+)} \biggr) \leq C_T,
\end{equation}
for some $ C_T $ depending only on $ T $, $ v^\mrm o, w^\mrm o, V$, and the initial data. Consequently, one has that 
\begin{equation}\label{est:007}
    \lim_{\eta \rightarrow \infty} (v, \theta) = 0.
\end{equation}
It is easy to verify that, from \eqref{eq:bl-app-1-v} and \eqref{eq:bl-app-1-w}, 
\begin{equation}
    % \begin{aligned}
    \nablah \pi  = - \int_0^1 \biggl(\dt v + (V + v^\mrm o)\cdot \nablah v + w^\mrm o \partial_\eta v + v \cdot \nablah V + \int_0^\eta \nablah \theta (s)\,ds\biggr) \,d\eta 
    % \end{aligned}
\end{equation}
is bounded in $ H^{s-1}(\mathbb T^2) $ for a.e. $ t \in (0,T) $. Therefore, one can send $ \eta \rightarrow \infty $ in \eqref{eq:bl-app-1-v}, and conclude that
\begin{equation}
    \int_0^\infty \nablah \theta (s) \,ds \in H^{s-3/2}(\mathbb T^2), \qquad a.e. ~ t \in (0,T).
\end{equation} 
After applying dominate convergence theorem to $ h_\eta(x,y):= \nablah \pi(x,y) + \int_0^\eta \nablah \theta(x,y,s)\,ds $ as $ \eta \rightarrow \infty $, one can obtain that
\begin{equation}\label{sol:pi}
    \nablah \pi  = - \int_0^\infty \nablah \theta(s)\,ds 
\end{equation}
and therefore \eqref{eq:bl-app-1-v} can be written as
\begin{equation*}\tag{\ref{eq:bl-app-1-v}'}\label{eq:bl-app-l-1-v-1}
    \dt v + (V + v^\mrm o)\cdot \nablah v + w^\mrm o \partial_\eta v + v \cdot \nablah V - \int_\eta^\infty \nablah \theta (s)\,ds = 0.
\end{equation*}

\bigskip
Now we have obtained an iteration scheme, i.e., \eqref{sys:bl-app-1}, from which we shall construct a contraction mapping from $ (v^\mrm o,\theta^\mrm o) $ to $ (v, \theta) $ in the next subsection. The fixed point obtained by such an iteration will be the unique solution to the boundary layer system \eqref{sys:bl}.

\subsection{The contraction mapping}\label{sec:mapping}

For the sake of clear presentation, we will assume $ V \equiv 0 $ in this section. The general case when $ V \neq 0 $ follows easily with proper modification. 

\subsubsection{Preliminaries and the main theorem}\label{sec:functional_space_bl}

As mentioned before, due to the loss of regularity in the horizontal directions, we shall consider solutions to system \eqref{sys:bl}, and system \eqref{sys:bl-app-1} in the analytic functional space. Our choice of the analytic norm should capture the following two characteristics of the systems:
\begin{enumerate}
    \item the boundedness of $ w $ and $ w^\mrm o $;
    \item the boundedness of $ \int_\eta^\infty \nablah \theta(s)\,ds $.
\end{enumerate}
A natural candidate is to consider $ (v,\theta) $ such that $ (v(\eta), \theta(\eta)) \simeq e^{-\beta \eta} $ for any $ \beta > 0 $, i.e., exponential decay in the $ \eta $-variable. Under this assumption, one will have that $ e^{\beta\eta} \int_\eta^\infty \nablah \theta(s)\,ds $ is uniformly bounded with respect to the $ \eta $-variable. 

Owing to the observation above, we consider the following analytic norm:
For any $ m \in \mathbb N $ and $ d > 0 $, denote by 
\begin{equation}\label{def:semi-norm}
    \vert g \vert_{d, m}:= \sum_{\vert\alpha\vert = m} \sup_{\eta \geq 0} e^{d\eta}\norm{\partial^\alpha g(\eta)}{L^2(\mathbb T^2)} = \sum_{\vert\alpha\vert = m} \norm{e^{d\eta} \partial^\alpha g}{L_\eta^\infty(\mathbb R^+, L_{x,y}^2(\mathbb T^2))}
\end{equation}
for multi-index $ \alpha := (\alpha_1,\alpha_2,\alpha_3) $. Here
\begin{equation}\label{def:multi-index-diff}
    \partial^\alpha g = \partial_x^{\alpha_1}\partial_y^{\alpha_2}\partial_\eta^{\alpha_3} g.
\end{equation}
Then for $ r \geq 0 $ and $ \tau > 0 $, we define the space of real-analytic functions as 
\begin{equation}\label{def:analytic-space}
    X_\tau := \biggl\lbrace v,\theta \in C^\infty(\Omega) :  \int_0^\infty \dvh v(\eta) \,d\eta = 0, ~ \norm{v,\theta}{X_\tau} < \infty \biggr\rbrace,
\end{equation}
where
\begin{equation}\label{def:analytic-norm}
    \norm{v,\theta}{X_\tau} := \sum_{m=0}^\infty \vert v, \theta \vert_{d, m} \dfrac{(m+1)^r \tau^m}{m!}.
\end{equation}
We have used and will use the notations $ \norm{v,\theta}{X} = \norm{v}{X} + \norm{\theta}{X} $ and $ \vert v, \theta \vert_X = \vert v \vert_X + \vert \theta \vert_X $ to denote semi-norms and norms of many functions. similarly, denote by
\begin{equation}\label{def:analytic-gaining-space}
    Y_\tau := \biggl\lbrace v, \theta \in X_\tau (\Omega):  \norm{v,\theta}{Y_\tau} < \infty \biggr\rbrace,
\end{equation}
where
\begin{equation}\label{def:analytic-gaining-norm}
    \norm{v,\theta}{Y_\tau} := \sum_{m=1}^\infty \vert v,\theta \vert_{d, m}\dfrac{(m+1)^r\tau^{m-1}}{(m-1)!}.
\end{equation}

\bigskip

With the notations above, we can now state the main theorem of this paper in full details as follows:
\begin{thm}\label{thm:b-l-detailed}
	For $ r \geq 2 $, there exists a unique solution $ (v, \theta) \in L^\infty(0,T;X_{\tau(t)})\cap L^1(0,T;Y_{\tau(t)}) $ to the boundary layer system \eqref{sys:bl}, with $ \tau $ determined by $ \eqref{eq:tau} $ and $ \tau(t=0) = \tau_0 > 0 $. Here $ T >0 $ depends on $ \tau_0 > 0 $ and $ \norm{v_\mrm{in},\theta_\mrm{in}}{X_{\tau_0}} $. See \eqref{def:mapping-space} for the precise estimates. 
\end{thm}
Here, $ r \geq 2 $ is to ensure the inequalities \eqref{ineq:0911-1}, \eqref{ineq:0911-2}, \eqref{ineq:0911-3}, and \eqref{ineq:0911-4}, below, hold. 

\subsubsection{Boundedness estimate}\label{sec:boundedness-estimates}

After applying $ \partial^\alpha $ to \eqref{eq:bl-app-l-1-v-1} and \eqref{eq:bl-app-1-theta} and taking the $ L^2 $-inner product in the horizontal variables $(x,y)$ with $ \partial^\alpha v $ and $ \partial^\alpha \theta $, respectively, one obtains
\begin{equation}\label{est:101}
    \begin{gathered}
    \dfrac{1}{2}\dfrac{d}{dt} \norm{\partial^\alpha v(\eta), \partial^\alpha \theta(\eta)}{L^2(\mathbb T^2)}^2 = - \int_{\mathbb T^2} \partial^\alpha ( v^\mrm o \cdot \nablah v + w^\mrm o \partial_\eta v ) \cdot \partial^\alpha v \,dxdy \\
    - \int_{\mathbb T^2} \partial^\alpha ( v^\mrm o \cdot \nablah \theta + w^\mrm o \partial_\eta \theta ) \times \partial^\alpha \theta \,dxdy\\
    + \int_{\mathbb T^2} \biggl( \partial^\alpha \int_\eta^\infty \nablah  \theta(s) \,ds \cdot \partial^\alpha v - \partial^\alpha w \times \partial^\alpha \theta \biggr) \,dxdy\\
    \leq \sum_{f = v,\theta} \norm{\partial^\alpha(v^\mrm o \cdot \nablah f + w^\mrm o \partial_\eta f)(\eta) }{L^2(\mathbb T^2)}
    \times \norm{\partial^\alpha v(\eta), \partial^\alpha \theta(\eta)}{L^2(\mathbb T^2)} \\
    + \norm{\partial^\alpha \int_\eta^\infty \nablah \theta(s)\,ds, \partial^\alpha w(\eta)}{L^2(\mathbb T^2)} \times \norm{\partial^\alpha v(\eta), \partial^\alpha \theta(\eta)}{L^2(\mathbb T^2)}.
    \end{gathered} 
\end{equation}
Therefore, one can write down, from \eqref{def:semi-norm}, \eqref{def:analytic-norm}, and \eqref{def:analytic-gaining-norm}, that
\begin{equation}\label{est:102}
    \begin{gathered}
    \dfrac{d}{dt} \norm{v, \theta}{X_\tau} \leq \dot\tau \norm{v,\theta}{Y_\tau} + \sum_{f = v,\theta} \biggl( \norm{v^\mrm o \cdot \nablah f}{X_\tau} + \norm{w^\mrm o \partial_\eta f}{X_\tau} \biggr) \\
     + \norm{\int_\eta^\infty \nablah \theta(s)\,ds, w}{X_\tau}.
    \end{gathered}
\end{equation}
Now we estimate the right hand side of \eqref{est:102}.

\bigskip

{\noindent\bf Estimate of $ \norm{v^\mrm o \cdot \nablah f}{X_\tau} $} Applying Leibniz's rule and triangle inequality of the $ L^2(\mathbb T^2) $-norm, one has that 
\begin{equation}
    \begin{aligned}
        & \norm{v^\mrm o \cdot \nablah f}{X_\tau} \leq \sum_{m=0}^\infty \sum_{\substack{\beta \leq \alpha\\\vert \alpha \vert = m}} \binom{\alpha}{\beta} \sup_{\eta \geq 0} e^{d\eta} \norm{\partial^\beta v^\mrm o\cdot \nablah \partial^{\alpha-\beta}f}{L^2(\mathbb T^2)} \dfrac{(m+1)^r\tau^m}{m!} \\
        & \qquad = \underbrace{\sum_{m=0}^\infty\sum_{j=0}^{[m/2]} \sum_{\substack{\beta \leq \alpha, |\beta| = j\\\vert \alpha \vert = m}} \binom{\alpha}{\beta} \sup_{\eta \geq 0} e^{d\eta} \norm{\partial^\beta v^\mrm o\cdot \nablah \partial^{\alpha-\beta}f}{L^2(\mathbb T^2)} \dfrac{(m+1)^r\tau^m}{m!}}_{=:\mathcal U_\mrm{low}} \\
        & \qquad + \underbrace{\sum_{m=0}^\infty\sum_{j=[m/2]+1}^{m} \sum_{\substack{\beta \leq \alpha, |\beta| = j\\\vert \alpha \vert = m}} \binom{\alpha}{\beta} \sup_{\eta \geq 0} e^{d\eta} \norm{\partial^\beta v^\mrm o\cdot \nablah \partial^{\alpha-\beta}f}{L^2(\mathbb T^2)} \dfrac{(m+1)^r\tau^m}{m!}}_{=:\mathcal U_\mrm{high}}.
    \end{aligned}
\end{equation}

In the case when $ 0 \leq |\beta| \leq |\alpha - \beta| $, applying the two-dimensional Agmon inequality leads to
\begin{equation}
    \begin{aligned}
        & \norm{\partial^\beta v^\mrm o\cdot \nablah \partial^{\alpha-\beta}f}{L^2(\mathbb T^2)} \leq C \norm{\partial^\beta v^\mrm o}{L^\infty(\mathbb T^2)} \norm{\nablah \partial^{\alpha-\beta}f}{L^2(\mathbb T^2)}\\
        & \qquad \leq C \norm{\partial^\beta v^\mrm o}{L^2(\mathbb T^2)}^{1/2}\norm{\Deltah\partial^\beta v^\mrm o}{L^2(\mathbb T^2)}^{1/2}\norm{\nablah \partial^{\alpha-\beta}f}{L^2(\mathbb T^2)} \\ 
        & \qquad\qquad + C \norm{\partial^{\beta}v^\mrm o}{L^2(\mathbb T^2)} \norm{\nablah \partial^{\alpha-\beta}f}{L^2(\mathbb T^2)}.
    \end{aligned}
\end{equation}
Notice that, for $ \beta \leq \alpha $, it holds $ |\alpha-\beta| = |\alpha| - |\beta| $, and 
\begin{equation}\label{ineq:101}
    \binom{\alpha}{\beta}\leq \binom{|\alpha|}{|\beta|} \quad \text{as well as} \quad \sum_{\substack{\beta\leq \alpha, |\beta|=j \\ |\alpha|=m}}a_\beta b_{\alpha-\beta} = \sum_{|\beta|=j}a_\beta \times \sum_{|\gamma| = m-j}b_\gamma.
\end{equation}
Therefore, one has
\begin{equation}\label{est:103}
    \begin{aligned}
        & \mathcal U_\mrm{low} \leq C \sum_{m=0}^\infty \sum_{j=0}^{[m/2]}\sum_{\substack{\beta \leq \alpha, |\beta| = j\\\vert \alpha \vert = m}} \sup_{\eta\geq 0}e^{d\eta} \biggl\lbrace \binom{\alpha}{\beta} \norm{\partial^\beta v^\mrm o}{L^2(\mathbb T^2)}^{1/2}\norm{\Deltah\partial^\beta v^\mrm o}{L^2(\mathbb T^2)}^{1/2}\\
        & \qquad\qquad\qquad \times \norm{\nablah \partial^{\alpha-\beta}f}{L^2(\mathbb T^2)} \dfrac{(m+1)^r\tau^m}{m!} \biggr\rbrace \\
        & \qquad\qquad + C \sum_{m=0}^\infty \sum_{j=0}^{[m/2]}\sum_{\substack{\beta \leq \alpha, |\beta| = j\\\vert \alpha \vert = m}} \sup_{\eta\geq 0}e^{d\eta} \biggl\lbrace \binom{\alpha}{\beta} \norm{\partial^{\beta}v^\mrm o}{L^2(\mathbb T^2)} \\
        & \qquad\qquad\qquad \times \norm{\nablah \partial^{\alpha-\beta}f}{L^2(\mathbb T^2)} \dfrac{(m+1)^r\tau^m}{m!} \biggr\rbrace \\
        & \leq C \sum_{m=0}^\infty \sum_{j=0}^{[m/2]} |v^\mrm o|^{1/2}_{d,j}|v^\mrm o|^{1/2}_{d,j+2}|f|_{d,m-j+1}
        %\sup_{\eta\geq 0}e^{d\eta} (\sum_{|\beta|=j}\norm{\partial^\beta v^\mrm o}{L^2(\mathbb T^2)})^{1/2} (\sum_{|\beta|=j}\norm{\Deltah\partial^\beta v^\mrm o}{L^2(\mathbb T^2)})^{1/2}\sum_{|\gamma| = m-j}\norm{\nablah\partial^\gamma f}{L^2(\mathbb T^2)}
        \binom{m}{j}\dfrac{(m+1)^r\tau^m}{m!} \\
        & \qquad + C \sum_{m=0}^\infty \sum_{j=0}^{[m/2]}|v^\mrm o|_{d,j}|f|_{d,m-j+1} \binom{m}{j}\dfrac{(m+1)^r\tau^m}{m!},
    \end{aligned}
\end{equation}
where we have applied, thanks to \eqref{ineq:101} and the Cauchy-Schwarz inequality,
\begin{equation}
    \begin{gathered}
    \sum_{\substack{\beta \leq \alpha, |\beta| = j\\\vert \alpha \vert = m}} \binom{\alpha}{\beta} \norm{\partial^\beta v^\mrm o}{L^2(\mathbb T^2)}^{1/2}\norm{\Deltah\partial^\beta v^\mrm o}{L^2(\mathbb T^2)}^{1/2}\norm{\nablah \partial^{\alpha-\beta}f}{L^2(\mathbb T^2)} \\
    \leq \binom{m}{j} (\sum_{|\beta|=j}\norm{\partial^\beta v^\mrm o}{L^2(\mathbb T^2)})^{1/2} (\sum_{|\beta|=j}\norm{\Deltah\partial^\beta v^\mrm o}{L^2(\mathbb T^2)})^{1/2}\sum_{|\gamma| = m-j}\norm{\nablah\partial^\gamma f}{L^2(\mathbb T^2)}, \\
    \text{and}\\
    \sum_{\substack{\beta \leq \alpha, |\beta| = j\\\vert \alpha \vert = m}} \binom{\alpha}{\beta} \norm{\partial^{\beta}v^\mrm o}{L^2(\mathbb T^2)} \norm{\nablah \partial^{\alpha-\beta}f}{L^2(\mathbb T^2)}\\
    \leq \binom{m}{j}\sum_{|\beta| = j}\norm{\partial^\beta v^\mrm o}{L^2(\mathbb T^2)} \times \sum_{|\beta|=m-j} \norm{\nablah\partial^\gamma f}{L^2(\mathbb T^2)}.
    \end{gathered}
\end{equation}
In addition, for $ m \geq 0, 0\leq j \leq [m/2], r \geq 1 $, one can verify the following inequality:
\begin{equation}\label{ineq:0911-1}
    \begin{gathered}
    \binom{m}{j} \dfrac{(m+1)^r}{m!}\dfrac{(m-j)!}{(m-j+2)^r} \dfrac{j!^{1/2}(j+2)!^{1/2}}{(j+1)^{r/2}(j+3)^{r/2}} \\ = \dfrac{(m+1)^r}{(m-j+2)^r}\dfrac{(j+1)^{1/2}(j+2)^{1/2}}{(j+1)^{r/2}(j+3)^{r/2}} \leq C,
    \end{gathered}
\end{equation}
for some constant $ C $ independent of $ m, j, r $. Thus $ \mathcal U_\mrm{low} $ from \eqref{est:103} can be further estimated as
\begin{equation}\label{est:104}
    \begin{aligned}
        & \mathcal U_\mrm{low} \leq C \sum_{m=0}^\infty \sum_{j=0}^{[m/2]} \biggl( |v^\mrm o|_{d,j}\dfrac{(j+1)^r \tau^j}{j!} \biggr)^{1/2} \biggl( |v^\mrm o|_{d,j+2} \dfrac{(j+3)^r \tau^{j+2}}{(j+2)!} \biggr)^{1/2} \\
        & \qquad\qquad \times \biggl(|f|_{d,m-j+1}\dfrac{(m-j+2)^r\tau^{m-j}}{(m-j)!} \biggr) \tau^{-1} \\
        & \qquad + C \sum_{m=0}^\infty \sum_{j=0}^{[m/2]} \biggl(| v^\mrm o|_{d,j} \dfrac{(j+1)^r \tau^j}{j!} \biggr) \biggl( |f|_{d,m-j+1} \dfrac{(m-j+2)^r\tau^{m-j}}{(m-j)!} \biggr)\\
        & \leq C (1+\tau^{-1}) \norm{v^\mrm o}{X_\tau} \norm{f}{Y_\tau},
    \end{aligned}
\end{equation}
where we apply the Cauchy-Schwarz inequality in the last inequality.

On the other hand, in the case when $ |\alpha-\beta| < |\beta| \leq |\alpha| $, one has
\begin{equation}
\begin{aligned}
    & \norm{\partial^\beta v^\mrm o\cdot \nablah \partial^{\alpha-\beta}f}{L^2(\mathbb T^2)} \leq C \norm{\partial^\beta v^\mrm o}{L^2(\mathbb T^2)} \norm{\nablah \partial^{\alpha-\beta}f}{L^\infty(\mathbb T^2)}\\
    & \qquad \leq C \norm{\partial^\beta v^\mrm o}{L^2(\mathbb T^2)}\norm{\nablah\partial^{\alpha-\beta} f}{L^2(\mathbb T^2)}^{1/2}\norm{\Deltah\nablah \partial^{\alpha-\beta}f}{L^2(\mathbb T^2)}^{1/2} \\ 
    & \qquad\qquad + C \norm{\partial^{\beta}v^\mrm o}{L^2(\mathbb T^2)} \norm{\nablah \partial^{\alpha-\beta}f}{L^2(\mathbb T^2)}.
\end{aligned}
\end{equation}
Then following similar arguments as above yields
\begin{equation}\label{est:105}
    \begin{aligned}
        & \mathcal U_{\mrm{high}} \leq C \sum_{m=0}^\infty \sum_{j=[m/2]+1}^m |v^\mrm o|_{d,j}|f|_{d,m-j+1}^{1/2}|f|_{d,m-j+3}^{1/2} \binom{m}{j}\dfrac{(m+1)^r\tau^m}{m!} \\
        & \qquad\qquad + C \sum_{m=0}^\infty \sum_{j=[m/2]+1}^m |v^\mrm o|_{d,j}|f|_{d,m-j+1}\binom{m}{j}\dfrac{(m+1)^r\tau^m}{m!}\\
        & \leq C(1+\tau^{-1})\norm{v^\mrm o}{X_\tau}\norm{f}{Y_\tau},
    \end{aligned}
\end{equation}
where the following inequality is applied:
\begin{equation}\label{ineq:0911-2}
\begin{gathered}
    \binom{m}{j}\dfrac{(m+1)^r}{m!}\dfrac{j!}{(j+1)^r}\dfrac{(m-j)!^{1/2}(m-j+2)!^{1/2}}{(m-j+2)^{r/2}(m-j+4)^{r/2}} \\
    = \dfrac{(m+1)^r}{(j+1)^r} \dfrac{(m-j+1)^{1/2}(m-j+2)^{1/2}}{(m-j+2)^{r/2}(m-j+4)^{r/2}} \leq C, 
\end{gathered}
\end{equation}
for $ m \geq 0, [m/2]+1 \leq j \leq m, r \geq 1 $. 

In summary, \eqref{est:104} and \eqref{est:105} yields
\begin{equation}\label{est:106}
    \sum_{f=v, \theta}\norm{v^\mrm o \cdot \nablah f}{X_\tau} \leq C (1+\tau^{-1}) \norm{v^\mrm o}{X_\tau}\norm{v, \theta}{Y_\tau}.
\end{equation}

\bigskip 
{\noindent\bf Estimate of $ \norm{w^\mrm o\partial_\eta f}{X_\tau} $} As before, applying Leibniz's rule and triangle inequality of the $ L^2(\mathbb T^2) $-norm, one has that
\begin{equation}
    \begin{aligned}
        & \norm{w^\mrm o\partial_\eta f}{X_\tau} \leq \sum_{m=0}^\infty \sum_{\substack{\beta \leq \alpha \\ |\alpha|=m}}\binom{\alpha}{\beta}\sup_{\eta \geq 0}e^{d\eta} \norm{w^\mrm o\partial_\eta f}{L^2(\mathbb T^2)}\dfrac{(m+1)^r\tau^m}{m!}\\
        & \qquad = \underbrace{\sum_{m=0}^\infty \sum_{j=0}^{[m/2]} \sum_{\substack{\beta \leq \alpha, |\beta| = j \\ |\alpha|=m}}\binom{\alpha}{\beta}\sup_{\eta \geq 0}e^{d\eta} \norm{w^\mrm o\partial_\eta f}{L^2(\mathbb T^2)}\dfrac{(m+1)^r\tau^m}{m!}}_{=: \mathcal V_\mrm{low}} \\
        & \qquad + \underbrace{\sum_{m=0}^\infty \sum_{j=[m/2]+1}^{m} \sum_{\substack{\beta \leq \alpha, |\beta| = j \\ |\alpha|=m}}\binom{\alpha}{\beta}\sup_{\eta \geq 0}e^{d\eta} \norm{w^\mrm o\partial_\eta f}{L^2(\mathbb T^2)}\dfrac{(m+1)^r\tau^m}{m!}}_{=: \mathcal V_\mrm{high}}.
    \end{aligned}
\end{equation}
To handle $ w^\mrm o $, we need to capture the decay-in-$ \eta $ of $ v^\mrm o $. Indeed, thanks to \eqref{est:005}, one can write
\begin{equation}
    w^\mrm o(\eta) = \int_0^\eta \partial_\eta w^\mrm o(s)\,ds = - \int_0^\eta \dvh v^\mrm o(s)\,ds = \int_\eta^\infty \dvh v^\mrm o(s)\,ds.
\end{equation}
Therefore, for $ \beta = (\beta_1, \beta_2, 0 ) $, one has
\begin{equation}\label{est:107}
    \norm{\partial^\beta w^\mrm o(\eta)}{L^2(\mathbb T^2)} \leq \dfrac{e^{-d\eta}}{d} \sup_{s\geq 0} e^{d s}\norm{\partial^\beta \dvh v^\mrm o(s)}{L^2(\mathbb T^2)},
\end{equation}
and for $ \beta = (\beta_1, \beta_2, \beta_3) $ with $ \beta_3 \geq 1 $, it follows 
\begin{equation}\label{est:108}
    \norm{\partial^\beta w^\mrm o(\eta)}{L^2(\mathbb T^2)} = % \norm{\partial^{\beta-(0,0,1)} \partial_\eta w^\mrm o(\eta)}{L^2(\mathbb T^2)} \leq 
    \norm{\partial^{\beta-(0,0,1)}\dvh v^\mrm o(\eta)}{L^2(\mathbb T^2)}.
\end{equation}
Thus repeating the previous arguments on $ \mathcal U_\mrm{low} $ implies that 
\begin{equation}\label{est:109}
    \begin{aligned}
        & \mathcal V_\mrm{low} \leq C \sum_{m=0}^\infty \sum_{j=0}^{[m/2]} \sum_{\substack{\beta \leq \alpha, |\beta|=j \\ |\alpha|=m}} \sup_{\eta \geq 0} e^{d\eta} \biggl\lbrace \binom{\alpha}{\beta} \norm{\partial^\beta w^\mrm o}{L^2(\mathbb T^2)}^{1/2}\norm{\Deltah\partial^\beta w^\mrm o}{L^2(\mathbb T^2)}^{1/2} \\
        & \qquad\qquad\qquad \times \norm{\partial_\eta \partial^{\alpha - \beta} f}{L^2(\mathbb T^2)} \dfrac{(m+1)^r \tau^m}{m!} \biggr\rbrace\\
        & \qquad\qquad + C \sum_{m=0}^\infty \sum_{j=0}^{[m/2]} \sum_{\substack{\beta \leq \alpha, |\beta|=j \\ |\alpha|=m}} \sup_{\eta \geq 0} e^{d\eta} \biggl\lbrace \binom{\alpha}{\beta} \norm{\partial^\beta w^\mrm o}{L^2(\mathbb T^2)} \\
        & \qquad\qquad\qquad \times \norm{\partial_\eta \partial^{\alpha-\beta}f}{L^2(\mathbb T^2)} \dfrac{(m+1)^r \tau^m}{m!} \biggr\rbrace \\
        & \qquad \leq C_d \sum_{m=0}^\infty \sum_{j=0}^{[m/2]} (|v^\mrm o|_{d,j}^{1/2}+|v^\mrm o|_{d,j+1}^{1/2})(|v^\mrm o|_{d,j+2}^{1/2}+|v^\mrm o|_{d,j+3}^{1/2}) |f|_{d,m-j+1} \\
        & \qquad\qquad\qquad \times \binom{m}{j}\dfrac{(m+1)^r\tau^m}{m!} \\
        & \qquad\qquad + C_d \sum_{m=0}^\infty \sum_{j=0}^{[m/2]} (|v^\mrm o|_{d,j}+|v^\mrm o|_{d,j+1}) |f|_{d,m-j+1} 
        % & \qquad\qquad\qquad \times 
        \binom{m}{j}\dfrac{(m+1)^r\tau^m}{m!},
    \end{aligned}
\end{equation}
where we have applied \eqref{est:107} and \eqref{est:108}. Furthermore, for $ m \geq 0, 0\leq j \leq [m/2], r \geq 2 $, the following inequality holds:
\begin{equation}\label{ineq:0911-3}
    \begin{gathered}
        \binom{m}{j}\dfrac{(m+1)^r}{m!} \dfrac{(m-j)!}{(m-j+2)^r} \dfrac{(j+1)!^{1/2}(j+3)!^{1/2}}{(j+2)^{r/2}(j+4)^{r/2}} \\
        = \dfrac{(m+1)^r}{(m-j+2)^r} \dfrac{(j+1)(j+2)^{1/2}(j+3)^{1/2}}{(j+2)^{r/2}(j+4)^{r/2}} \leq C. 
    \end{gathered}
\end{equation}
Therefore, one can obtain from \eqref{est:109}, \eqref{def:analytic-norm}, and \eqref{def:analytic-gaining-norm}, that
\begin{equation}\label{est:110}
    \mathcal V_\mrm{low} \leq C_d(1+\tau^{-2}) \norm{v^\mrm o}{X_\tau} \norm{f}{Y_\tau}.
\end{equation}
Meanwhile, one has
\begin{equation}\label{est:111}
    \begin{aligned}
        & \mathcal V_\mrm{high} \leq C_d \sum_{m=0}^\infty \sum_{j=[m/2]+1}^m (|v^\mrm o|_{d,j} + |v^\mrm o|_{d,j+1}) |f|^{1/2}_{d,m-j+1}|f|^{1/2}_{d,m-j+3}\\
        & \qquad\qquad\qquad \times \binom{m}{j} \dfrac{(m+1)^r\tau^m}{m!} \\
        & \qquad\qquad + C_d \sum_{m=0}^\infty \sum_{j=[m/2]+1}^m (|v^\mrm o|_{d,j} + |v^\mrm o|_{d,j+1}) |f|_{d,m-j+1}\\
        & \qquad\qquad\qquad \times \binom{m}{j} \dfrac{(m+1)^r\tau^m}{m!}.
    \end{aligned}
\end{equation}
The one containing $ |v^\mrm o|_{d,j} $ on the right hand side of \eqref{est:111} can be estimated by $ C_d (1+\tau^{-1})\norm{v^\mrm o}{X_\tau} \norm{f}{Y_\tau} $, similar to $ \mathcal U_\mrm{high} $. The one containing $ |v^\mrm o|_{d,j} $, on the other hand, can be estimated by $ C_d (\tau^{-1} + \tau^{-2}) \norm{v^\mrm o}{Y_\tau}\norm{f}{X_\tau} $ thanks to the following inequality:
\begin{equation}\label{ineq:0911-4}
    \begin{gathered}
        \binom{m}{j} \dfrac{(m+1)^r}{m!} \dfrac{j!}{(j+2)^r}\dfrac{(m-j+1)!^{1/2}(m-j+3)!^{1/2}}{(m-j+2)^{r/2}(m-j+4)^{r/2}}\\
        = \dfrac{(m+1)^r}{(j+2)^r} \dfrac{(m-j+1)(m-j+2)^{1/2}(m-j+3)^{1/2}}{(m-j+2)^{r/2}(m-j+4)^{r/2}} \leq C
    \end{gathered}
\end{equation}
for $ m \geq 0, [m/2]+1\leq j \leq m, r\geq 2 $. 

In summary, from \eqref{est:110} and \eqref{est:111}, we have shown 
\begin{equation}\label{est:112}
    \sum_{f=v,\theta} \norm{w^\mrm o \partial_\eta f}{X_\tau} \leq C_d(1+\tau^{-1}) \norm{v^\mrm o}{X_\tau} \norm{v, \theta}{Y_\tau} + C_d(\tau^{-1} + \tau^{-2}) \norm{v^\mrm o}{Y_\tau} \norm{v, \theta}{X_\tau}.
\end{equation}

\bigskip
{\noindent\bf Estimate of $ \norm{\int_\eta^\infty \nablah \theta(s)\,ds, w}{X_\tau} $} Similar to the estimates \eqref{est:107} and \eqref{est:108}, one can obtain
\begin{equation}
    \begin{aligned}
        &\sup_{\eta \geq 0}e^{d\eta}\norm{\partial^\beta \int_\eta^\infty \nablah \theta(s)\,ds, \partial^\beta w(\eta)}{L^2(\mathbb T^2)} \\
        & \leq \begin{cases}
            d^{-1} \sup_{s\geq 0} e^{ds}\norm{\partial^\beta \dvh v(s), \partial^\beta \nablah \theta(s)}{L^2(\mathbb T^2)} & \text{if $\beta_3 = 0$} \\
            \sup_{s\geq 0} e^{ds}\norm{\partial^{\beta-(0,0,1)} \dvh v(s), \partial^{\beta-(0,0,1)} \nablah \theta(s)}{L^2(\mathbb T^2)} & \text{if $\beta_3 \geq 1$}.
        \end{cases}
    \end{aligned}
\end{equation}
Therefore, we have that
\begin{equation}\label{est:113}
    \norm{\int_\eta^\infty \nablah \theta(s)\,ds, w}{X_\tau} \leq \norm{v,\theta}{X_\tau} + d^{-1} \norm{v,\theta}{Y_\tau}.
\end{equation}

\bigskip
{\noindent\bf Summary} From \eqref{est:106}, \eqref{est:112}, and \eqref{est:113}, one can conclude from \eqref{est:102} that 
\begin{equation}\label{est:114}
    \begin{gathered}
        \dfrac{d}{dt}\norm{v,\theta}{X_\tau} \leq \lbrack\dot \tau + C_d (1 + \tau^{-1})\norm{v^\mrm o}{X_\tau} + d^{-1} \rbrack \norm{v,\theta}{Y_\tau} \\
        + \lbrack 1 + C_d(\tau^{-1}+\tau^{-2})\norm{v^\mrm o}{Y_\tau}\rbrack \norm{v,\theta}{X_\tau}.
    \end{gathered}
\end{equation}

Consider $ M := \norm{v_\mrm{in},\theta_\mrm{in}}{X_{\tau_0}} $, where
\begin{equation}
    (v,\theta)\vert_{t=0} = (v_\mrm{in},\theta_\mrm{in}) \in X_{\tau_0}
\end{equation}
with $ \tau_0 > 0 $.
Let $ \tau $ satisfy $ \tau(t=0) = \tau_0 $ and 
\begin{equation}\label{eq:tau}
    \dot \tau (t) + 2 C_d(1+\tau^{-1}) M + d^{-1} + 4C_d(2\tau_0^{-1} + 4 \tau_0^{-2})M = - 4C_d(1+2\tau_0^{-1})M.
\end{equation}
Then $ \tau $ decreases in $ t $, and for $ T < 1/4 $ small enough, depending on $ M $ and $ \tau_0 $, we have that, for $0\leq t \leq T$,
\begin{equation}
    \tau_0 / 2 \leq  \tau(t) \leq \tau_0.
\end{equation}
Then integrating \eqref{est:114} implies 
\begin{equation}\label{est:115}
    \begin{gathered}
        \sup_{0\leq t\leq T} \norm{v(t),\theta(t)}{X_{\tau(t)}} + \biggl(4 C_d(2\tau_0^{-1} + 4\tau_0^{-2})M+4C_d(1+2\tau_0^{-1})M \biggr) \\
        \qquad\qquad\qquad \times \int_0^T \norm{v(t),\theta(t)}{Y_{\tau(t)}} \,dt \\
        \leq M + C_d(1+2\tau_0^{-1}) (\sup_{0\leq t \leq T}\norm{v^\mrm o(t)}{X_{\tau(t)}} - 2M)\int_0^T \norm{v(t),\theta(t)}{Y_{\tau(t)}} \,dt  \\
        + \biggl(\dfrac{1}{4} + C_d(2\tau_0^{-1} + 4 \tau_0^{-2}) \int_0^T \norm{v^\mrm o(t)}{Y_{\tau(t)}}\,dt \biggr) \sup_{0\leq t \leq T} \norm{v(t),\theta(t)}{X_{\tau(t)}}.
    \end{gathered}
\end{equation}
Therefore, with $ \tau $ givens by \eqref{eq:tau} and $ T $ as above, the following mapping
\begin{equation}\label{def:mapping}
    \mathfrak M: (v^\mrm o,\theta^\mrm o) \mapsto (v, \theta)
\end{equation}
is bounded in 
\begin{equation}\label{def:mapping-space}
    \begin{gathered}
        X_{M,T}:= \biggl\lbrace v,\theta \in L^\infty(0,T;X_\tau)\cap L^1(0,T;Y_\tau): \sup_{0\leq t \leq T} \norm{v(t),\theta(t)}{X_{\tau(t)}} \\
        + (8 C_d(2\tau_0^{-1} + 4\tau_0^{-2})M+8C_d(1+2\tau_0^{-1})M) \int_0^T \norm{v(t),\theta(t)}{Y_{\tau(t)}}\,dt \leq 2 M, \\ 
        (v,\theta)\vert_{t=0} = (v_\mrm{in},\theta_\mrm{in}), ~ \norm{v_\mrm{in},\theta_\mrm{in}}{X_{\tau_0}} = M, ~ \int_0^\infty \dvh v(s)\,ds = 0
        \biggr\rbrace.
    \end{gathered}
\end{equation}

\subsubsection{Contracting and the fixed point}\label{sec:contracting}
Thanks to the boundedness estimate from section \ref{sec:boundedness-estimates}, to obtain a fixed point of $ \mathcal M $ in \eqref{def:mapping}, it suffices to show that $ \mathcal M $ is contracting in $L^\infty(0,T;X_\tau)\cap L^1(0,T;Y_\tau)$. Let $ (v_\mrm i,\theta_\mrm i) = \mathcal M(v_\mrm i^\mrm o, \theta_\mrm i^\mrm o) $, $ i = 1,2 $, and $ \delta \psi = \psi_1 - \psi_2 $ for $ \psi \in \lbrace v,\theta, v^\mrm o, \theta^\mrm o \rbrace $. Then one has, from \eqref{sys:bl-app-1}
\begin{subequations}
    \begin{gather}
        \begin{gathered}
            \dt \delta v + v^\mrm o_1 \cdot \nablah \delta v +  w^\mrm o_1 \partial_\eta \delta v + \delta v^\mrm o \cdot \nablah v_2 + \delta  w^\mrm o \partial_\eta v_2 \\
            - \int_\eta^\infty \nablah \delta \theta(s)\,ds = 0, 
        \end{gathered}\label{eq:ctt-1} \\
        \dt \delta \theta + v^\mrm o_1 \cdot \nablah \delta \theta + w^\mrm o_1 \partial_\eta \delta \theta + \delta v^\mrm o \cdot \nablah \theta_2 + \delta w^\mrm o \partial_\eta \theta_2 + \delta w = 0. \label{eq:ctt-2}
    \end{gather}
\end{subequations}
Here $ w_\mrm i, w_\mrm i^\mrm o, \delta w, \delta w^\mrm o $ are defined similar as $ w, w^\mrm o $ using the boundary condition and the incompressible condition as before. Then repeating the same arguments as in section \ref{sec:boundedness-estimates} yields
\begin{equation}\label{est:116}
    \begin{aligned}
        & \dfrac{d}{dt}\norm{\delta v, \delta \theta}{X_\tau} \leq \lbrack\dot\tau + C_d(1+\tau^{-1}) \norm{v^\mrm o_1}{X_\tau} + d^{-1}\rbrack \norm{\delta v, \delta \theta}{Y_\tau}  \\
        & \qquad + \lbrack 1 + C_d (\tau^{-1} + \tau^{-2}) \norm{v^\mrm o_1}{Y_\tau} \rbrack \norm{\delta v, \delta \theta}{X_\tau} \\
        & \qquad + C_d(1+\tau^{-1}) \norm{\delta v^\mrm o}{X_\tau}\norm{v_2,\theta_2}{Y_\tau}\\ 
        & \qquad + C_d (\tau^{-1} + \tau^{-2}) \norm{\delta v^\mrm o}{Y_\tau}\norm{v_2,\theta_2}{X_\tau}.
    \end{aligned}
\end{equation}
Then integrating in time of \eqref{est:116} implies
\begin{equation}\label{est:117}
    \begin{aligned}
        & \sup_{0\leq t \leq T}\norm{\delta v(t), \delta \theta(t)}{X_{\tau(t)}} + \biggl(4C_d (2\tau_0^{-1} + 4 \tau_0^{-2})M + 4 C_d (1+2\tau_0^{-1})M\biggr)\\
        & \qquad\qquad\qquad \times \int_0^T \norm{\delta v(t),\delta \theta(t)}{Y_{\tau(t)}}\,dt \\
        & \leq \biggl( \dfrac{1}{4} + C_d(2\tau_0^{-1} + 4\tau^{-2})\int_0^T \norm{v_1^\mrm o(t)}{Y_{\tau(t)}}\,dt \biggr) \sup_{0\leq t\leq T} \norm{\delta v(t), \delta \theta(t)}{X_{\tau(t)}} \\
        & \qquad + C_d(1+2\tau_0^{-1}) \int_0^T\norm{v_2(t),\theta_2(t)}{Y_{\tau(t)}}\,dt \sup_{0\leq t\leq T} \norm{\delta v^\mrm o(t)}{X_{\tau(t)}} \\
        & \qquad + C_d(2\tau_0^{-1} + 4\tau_0^{-2}) \sup_{0\leq t \leq T} \norm{v_2(t),\theta_2(t)}{X_{\tau(t)}} \int_0^T \norm{\delta v^\mrm o}{Y_{\tau(t)}} \,dt
    \end{aligned}
\end{equation}
and therefore, substituting the boundedness estimate in \eqref{def:mapping-space} yield
\begin{equation}
    \begin{gathered}
    \dfrac{1}{2}\sup_{0\leq t \leq T}\norm{\delta v(t),\delta\theta(t)}{X_{\tau(t)}} + \biggl(4C_d (2\tau_0^{-1} + 4 \tau_0^{-2})M + 4 C_d (1+2\tau_0^{-1})M\biggr)\\
    \qquad\qquad\qquad\times \int_0^T \norm{\delta v(t),\delta \theta(t)}{Y_{\tau(t)}}\,dt\\
    \leq \dfrac{1}{4} \sup_{0\leq t \leq T}\norm{\delta v^\mrm o(t)}{X_{\tau(t)}} + 2C_d(2\tau_0^{-1} + 4 \tau_0^{-2})M \int_0^T \norm{\delta v^\mrm o(t),\delta \theta^\mrm o(t)}{Y_{\tau(t)}}\,dt.
\end{gathered}
\end{equation}
Consequently, $ \mathfrak M $ in \eqref{def:mapping} is contracting with respect to $ \mathfrak E $ defined by 
\begin{equation}\label{def:contracting-norm}
    \begin{gathered}
        \mathfrak E(v,\theta):= \dfrac{1}{2} \sup_{0\leq t \leq T}\norm{v(t),\theta(t)}{X_{\tau(t)}}  \\
        + 4C_d(2\tau_0^{-1} + 4 \tau_0^{-2})M \int_0^T \norm{v(t),\theta(t)}{Y_{\tau(t)}}\,dt,
    \end{gathered}
\end{equation}
and there is a fixed point of $ \mathfrak M $, which is the solution in $ X_{M,T} $ (defined in \eqref{def:mapping-space}) to the following system
\begin{subequations}
    \begin{gather}
        \dt v + v \cdot \nablah v + w \partial_\eta v - \int_\eta^\infty \nablah \theta(s)\,ds = 0, \\
        \dt \theta + v \cdot \nablah \theta + w \partial_\eta \theta + w = 0, \\
        \dvh v + \partial_\eta w = 0,
    \end{gather}
    with 
    \begin{equation}
        w\vert_{\eta = 0} \qquad \text{and} \qquad \int_0^\infty \dvh v(\eta)\,d\eta = 0.
    \end{equation}
\end{subequations}
This finishes the proof of theorem \ref{thm:b-l-detailed} for $ V \equiv 0 $. For $ V \in L^\infty(0,T;X_{\tau_0/2}) \cap L^1(0,T; Y_{\tau_0/2}) $, repeating the above arguments yields the solution to \eqref{sys:bl}. This finishes the proof of theorem \ref{thm:main-theorem-informal}.

\section*{Acknowledgement}

R.K.\ thanks Deutsche Forschungsgemeinschaft for funding under grant
CRC 1114 ``Scaling Cascades in Complex Systems'', Project No.\ 235221301, Project (C06) ``Multi-scale structure of atmospheric vortices'', and under grant FOR 5528 ``Mathematical Study of Geophysical Flow Models: Analysis and Computation'', Project No.\ 500072446, Project TP2 ``Scale Analysis and Asymptotic Reduced Models for the Atmosphere''.

%\pagebreak

\bibliographystyle{plain}
%\bibliography{ref.bib}

\end{document}